arXiv:1606.02118v3 [math.OC] 27 Oct 2016

# A Multi-step Inertial Forward–Backward Splitting Method for Non-convex Optimization

Jingwei Liang* Jalal M. Fadili* Gabriel Peyré†

**Abstract**

In this paper, we propose a multi-step inertial Forward–Backward splitting algorithm for minimizing the sum of two non-necessarily convex functions, one of which is proper lower semi-continuous while the other is differentiable with a Lipschitz continuous gradient. We first prove global convergence of the scheme with the help of the Kurdyka–Łojasiewicz property. Then, when the non-smooth part is also partly smooth relative to a smooth submanifold, we establish finite identification of the latter and provide sharp local linear convergence analysis. The proposed method is illustrated on a few problems arising from statistics and machine learning.

# 1 Introduction

## 1.1 Non-convex non-smooth optimization

Non-smooth optimization has proved extremely useful to all quantitative disciplines of science including statistics and machine learning. A common trend in modern science is the increase in size of datasets, which drives the need for more efficient optimization schemes. For large-scale problems with non-smooth and possibly non-convex terms, it is possible to generalize gradient descent with the Forward–Backward (FB) splitting scheme [4] (a.k.a proximal gradient descent), which includes projected gradient descent as a sub-case.

Formally, we equip $\mathbb{R}^n$ the $n$-dimensional Euclidean space with the standard inner product $\langle \cdot, \cdot \rangle$ and associated norm $\|\cdot\|$ respectively. Our goal is the generic minimization of composite objectives of the form

$$\min_{x \in \mathbb{R}^n} \big\{ \Phi(x) \stackrel{\text{def}}{=} R(x) + F(x) \big\}, \tag{$\mathcal{P}$}$$

where we have

(**A.1**) $R : \mathbb{R}^n \to \mathbb{R} \cup \{+\infty\}$ is the *penalty function* which is proper lower semi-continuous (lsc), and bounded from below;

(**A.2**) $F : \mathbb{R}^n \to \mathbb{R}$ is the *loss function* which is finite-valued, differentiable and its gradient $\nabla F$ is $L$-Lipschitz continuous.

Throughout, no convexity is imposed neither on $R$ nor on $F$.

The class of problems we consider is that of non-smooth and non-convex optimization problems. Here are some examples that are of particular relevance to problems in regression, machine learning and classification.

*Normandie University, ENSICAEN, UNICAEN, GREYC, E-mail: {Jingwei.Liang, Jalal.Fadili}@ensicaen.fr.
†CNRS, Ceremade, Université Paris-Dauphine, E-mail: Gabriel.Peyre@ceremade.dauphine.fr.

**Example 1.1 (Sparse regression).** Let $A \in \mathbb{R}^{m\times n}$, $y \in \mathbb{R}^m$, $\mu > 0$, and $\|x\|_0$ is the $\ell_0$ pseudo-norm (see Example 4.1). Consider (see *e.g.* [14])

$$\min_{x\in\mathbb{R}^n} \frac{1}{2}\|y - Ax\|^2 + \mu\|x\|_0. \tag{1.1}$$

**Example 1.2 (Principal component pursuit (PCP)).** The PCP problem [10] aims at decomposing a given matrix into *sparse* and *low-rank* components

$$\min_{(x_s,x_l)\in(\mathbb{R}^{n_1\times n_2})^2} \frac{1}{2}\|y - x_s - x_l\|_{\mathrm{F}}^2 + \mu_1\|x_s\|_0 + \mu_2 \mathrm{rank}(x_l), \tag{1.2}$$

where $\|\cdot\|_{\mathrm{F}}$ is the Frobenius norm and $\mu_1$ and $\mu_2 > 0$.

**Example 1.3 (Sparse Support Vector Machines).** One would like to find a linear decision function which minimizes the objective

$$\min_{(b,x)\in\mathbb{R}\times\mathbb{R}^n} \frac{1}{m}\sum_{i=1}^m G(\langle x,\, z_i\rangle + b, y_i) + \mu\|x\|_0, \tag{1.3}$$

where for $i = 1, \cdots, m$, $(z_i, y_i) \in \mathbb{R}^n \times \{\pm 1\}$ is the training set, and $G$ is a smooth loss function with Lipschitz-continuous gradient such as the squared hinge loss $G(\hat{y}_i, y_i) = \max(0, 1 - \hat{y}_i y_i)^2$ or the logistic loss $G(\hat{y}_i, y_i) = \log(1 + e^{-\hat{y}_i y_i})$.

**(Inertial) Forward–Backward** The Forward–Backward splitting method for solving ($\mathcal{P}$) reads

$$x_{k+1} \in \mathrm{prox}_{\gamma_k R}\big(x_k - \gamma_k \nabla F(x_k)\big), \tag{1.4}$$

where $\gamma_k > 0$ is a descent step-size, and

$$\mathrm{prox}_{\gamma R}(\cdot) \stackrel{\text{def}}{=} \mathrm{Argmin}_{x\in\mathbb{R}^n} \frac{1}{2}\|x - \cdot\|^2 + \gamma R(x), \tag{1.5}$$

denotes the proximity operator of $R$. $\mathrm{prox}_{\gamma R}(x)$ is non-empty under (**A.1**) and is set-valued in general. Lower-boundedness of $R$ can be relaxed by requiring e.g. coercivity of the objective in (1.5).

Since the pioneering work of Polyak [28] on the *heavy-ball method* approach to gradient descent, several works have adapted this methodology to various optimization schemes. For instance, the inertial proximal point algorithm [2, 3], or the inertial FB methods [26, 24, 5, 23]. The FISTA scheme [6, 11] also belongs to this class. See [23] for a detailed account.

**The non-convex case** In the context of non-convex optimization, [4] was the first to establish convergence of the FB iterates when the objective $\Phi$ satisfies the Kurdyka–Łojasiewicz property[1]. Following their footprints, [9, 27] established convergence of the special inertial schemes in [26] in the non-convex setting.

## 1.2 Contributions

In this paper, we introduce a novel inertial scheme (Algorithm 1) and study its global and local properties to solve the non-smooth and non-convex optimization problem ($\mathcal{P}$). More precisely, our main contributions can be summarized as follows.

[1]We are aware of the works existing on convergence of the objective sequence $\Phi(x_k)$ of FB, including rates, in the non-smooth and non-convex setting. But given that, in general, this does not say anything about convergence of the sequence of iterates $x_k$, they are irrelevant to our discussion.

**A globally convergent general inertial scheme** We propose a general multi-step inertial FB (MiFB) algorithm to solve ($\mathcal{P}$). This algorithm is very flexible as it allows higher memory and even *negative* inertial parameters (unlike previous work [23]). Global convergence of any bounded sequence of iterates to a critical point is proved when the objective $\Phi$ is lower-bounded and satisfies the Kurdyka–Łojasiewicz property.

**Local convergence properties under partial smoothness** Under the additional assumptions that the smooth part is locally $C^2$ around a critical point $x^\star$ (where $x_k \to x^\star$), and that the non-smooth component $R$ is partly smooth (see Definition 3.1) relative to an active submanifold $\mathcal{M}_{x^\star}$, we show that $\mathcal{M}_{x^\star}$ can be identified in finite time, *i.e.* $x_k \in \mathcal{M}_{x^\star}$ for all $k$ large enough. Building on this finite identification result, we provide a sharp local linear convergence analysis and we characterize precisely the corresponding convergence rate which, in particular, reveals the role of $\mathcal{M}_{x^\star}$. Moreover, this local convergence analysis naturally opens the door to higher-order acceleration, since under some circumstances, the original problem ($\mathcal{P}$) is eventually equivalent to locally minimizing $\Phi$ on $\mathcal{M}_{x^\star}$, and partial smoothness implies that $\Phi$ is actually $C^2$ on $\mathcal{M}_{x^\star}$.

---

**Algorithm 1:** A Multi-step Inertial Forward–Backward (MiFB)

---

**Initial**: $s \geq 1$ is an integer, $I = \{0, 1, \ldots, s-1\}$, $x_0 \in \mathbb{R}^n$ and $x_{-s} = \ldots = x_{-1} = x_0$.

**repeat**

Let $0 < \underline{\gamma} \leq \gamma_k \leq \bar{\gamma} < \frac{1}{L}$, $\{a_{0,k}, a_{1,k}, \ldots\} \in ]-1, 2]^s$, $\{b_{0,k}, b_{1,k}, \ldots\} \in ]-1, 2]^s$:

$$\begin{aligned} y_{a,k} &= x_k + \textstyle\sum_{i \in I} a_{i,k}(x_{k-i} - x_{k-i-1}), \\ y_{b,k} &= x_k + \textstyle\sum_{i \in I} b_{i,k}(x_{k-i} - x_{k-i-1}), \end{aligned} \tag{1.6}$$

$$x_{k+1} \in \mathrm{prox}_{\gamma_k R}\big(y_{a,k} - \gamma_k \nabla F(y_{b,k})\big). \tag{1.7}$$

$k = k + 1$;

**until** *convergence*;

---

## 1.3 Notations

Throughout the paper, $\mathbb{N}$ is the set of non-negative integers. For a nonempty closed convex set $\Omega \subset \mathbb{R}^n$, $\mathrm{ri}(\Omega)$ is its relative interior, and $\mathrm{par}(\Omega) = \mathbb{R}(\Omega - \Omega)$ is the subspace parallel to it.

Let $R : \mathbb{R}^n \to \mathbb{R} \cup \{+\infty\}$ be a lsc function, its domain is defined as $\mathrm{dom}(R) \stackrel{\mathrm{def}}{=} \{x \in \mathbb{R}^n : R(x) < +\infty\}$, and it is said to be proper if $\mathrm{dom}(R) \neq \emptyset$. We need the following notions from variational analysis, see e.g. [30] for details. Given $x \in \mathrm{dom}(R)$, the Fréchet subdifferential $\partial^F R(x)$ of $R$ at $x$, is the set of vectors $v \in \mathbb{R}^n$ that satisfies $\liminf_{z \to x,\, z \neq x} \frac{1}{\|x-z\|}(R(z) - R(x) - \langle v,\, z - x\rangle) \geq 0$. If $x \notin \mathrm{dom}(R)$, then $\partial^F R(x) = \emptyset$. The limiting-subdifferential (or simply subdifferential) of $R$ at $x$, written as $\partial R(x)$, is defined as $\partial R(x) \stackrel{\mathrm{def}}{=} \{v \in \mathbb{R}^n : \exists x_k \to x, R(x_k) \to R(x), v_k \in \partial^F R(x_k) \to v\}$. Denote $\mathrm{dom}(\partial R) \stackrel{\mathrm{def}}{=} \{x \in \mathbb{R}^n : \partial R(x) \neq \emptyset\}$. Both $\partial^F R(x)$ and $\partial R(x)$ are closed, with $\partial^F R(x)$ convex and $\partial^F R(x) \subset \partial R(x)$ [30, Proposition 8.5]. Since $R$ is lsc, it is (subdifferentially) regular at $x$ if and only if $\partial^F R(x) = \partial R(x)$ [30, Corollary 8.11].

An lsc function $R$ is $r$-prox-regular at $\bar{x} \in \mathrm{dom}(R)$ for $\bar{v} \in \partial R(\bar{x})$ if $\exists r > 0$ such that $R(x') > R(x) + \langle v,\, x' - x\rangle - \frac{1}{2r}\|x - x'\|^2$ $\forall x, x'$ near $\bar{x}$, $R(x)$ near $R(\bar{x})$ and $v \in \partial R(x)$ near $\bar{v}$.

A necessary condition for $x$ to be a minimizer of $R$ is $0 \in \partial R(x)$. The set of critical points of $R$ is $\mathrm{crit}(R) = \{x \in \mathbb{R}^n : 0 \in \partial R(x)\}$.

## 2 Global convergence of MiFB

This section is dedicated to the global convergence of the sequence generated by Algorithm 1.

**Kurdyka–Łojasiewicz property** Let $J : \mathbb{R}^n \to \mathbb{R} \cup \{+\infty\}$ be a proper lsc function. For $\eta_1, \eta_2$ such that $-\infty < \eta_1 < \eta_2 < +\infty$, define the set $[\eta_1 < J < \eta_2] \stackrel{\text{def}}{=} \{x \in \mathbb{R}^n : \eta_1 < J(x) < \eta_2\}$.

**Definition 2.1.** $J$ is said to have the Kurdyka–Łojasiewicz property at $\bar{x} \in \mathrm{dom}(J)$ if there exists $\eta \in ]0, +\infty]$, a neighbourhood $U$ of $\bar{x}$ and a continuous concave function $\varphi : [0, \eta[ \to \mathbb{R}_+$ such that

(i) $\varphi(0) = 0$, $\varphi$ is $C^1$ on $]0, \eta[$, and for all $s \in ]0, \eta[$, $\varphi'(s) > 0$;

(ii) for all $x \in U \cap [J(\bar{x}) < J < J(\bar{x}) + \eta]$, the Kurdyka–Łojasiewicz inequality holds

$$\varphi'\big(J(x) - J(\bar{x})\big)\mathrm{dist}\big(0, \partial J(x)\big) \geq 1. \tag{2.1}$$

Proper lsc functions which satisfy the Kurdyka–Łojasiewicz property at each point of $\mathrm{dom}(\partial J)$ are called KL functions.

Roughly speaking, KL functions become sharp up to reparameterization via $\varphi$, called a desingularizing function for $J$. Typical KL functions are the class of semi-algebraic functions, see [7, 8]. For instance, the $\ell_0$ pseudo-norm and the rank function (see Example 1.1-1.3 and Section 4.1) are indeed KL.

### 2.1 Global convergence

Let $\mu, \nu > 0$ be two constants. For $i \in I$ and $k \in \mathbb{N}$, define the following quantities,

$$\beta_k \stackrel{\text{def}}{=} \frac{1 - \gamma_k L - \mu - \nu\gamma_k}{2\gamma_k}, \ \underline{\beta} \stackrel{\text{def}}{=} \liminf_{k\in\mathbb{N}} \beta_k \ \text{and} \ \alpha_{k,i} \stackrel{\text{def}}{=} \frac{s a_{i,k}^2}{2\gamma_k \mu} + \frac{s b_{i,k}^2 L^2}{2\nu}, \ \bar{\alpha}_i \stackrel{\text{def}}{=} \limsup_{k\in\mathbb{N}} \alpha_{k,i}. \tag{2.2}$$

**Theorem 2.2 (Convergence of MiFB (Algorithm 1)).** *For problem ($\mathcal{P}$), suppose that* (**A.1**)-(**A.2**) *hold. Assume moreover that $\Phi$ is a proper lsc KL function which is bounded from below. For Algorithm 1, choose $\mu, \nu, \gamma_k, a_{i,k}, b_{i,k}$ such that*

$$\delta \stackrel{\text{def}}{=} \underline{\beta} - \textstyle\sum_{i\in I} \bar{\alpha}_i > 0. \tag{2.3}$$

*Then each bounded sequence $\{x_k\}_{k\in\mathbb{N}}$ satisfies*

(i) *$\{x_k\}_{k\in\mathbb{N}}$ has finite length, i.e. $\sum_{k\in\mathbb{N}} \|x_k - x_{k-1}\| < +\infty$;*

(ii) *There exists a critical point $x^\star \in \mathrm{crit}(\Phi)$ such that $\lim_{k\to\infty} x_k = x^\star$.*

(iii) *If $\Phi$ has the KL property at a global minimizer $x^\star$, then starting sufficiently close from $x^\star$, any sequence $\{x_k\}_{k\in\mathbb{N}}$ converges to a global minimum of $\Phi$ and satisfies (i).*

See the Section A for the detailed proof.

**Remark 2.3.**

(i) Boundedness of the sequence is automatically ensured under standard assumptions such as coercivity of $\Phi$.

(ii) Unlike existing work, *negative* inertial parameters are allowed by Theorem 2.2.

(iii) When $a_{i,k} \equiv 0$ and $b_{i,k} \equiv 0$, *i.e.* the case of FB splitting, condition (2.3) holds naturally as long as $\bar{\gamma} < \frac{1}{L}$ which recovers the case of [4];

(iv) From (2.2) and (2.3), we conclude the following:

(a) $s = 1$: if $b_{k,0} \equiv b, a_{k,0} \equiv a$ (*i.e.* constant inertial parameters), then (2.3) implies that $a, b$ must belong to an ellipsoid,

$$\frac{a^2}{2\underline{\gamma}\mu} + \frac{b^2}{2\nu/L^2} < \underline{\beta} = \frac{1-\bar{\gamma}L-\mu-\nu\bar{\gamma}}{2\bar{\gamma}}.$$

(b) When $s \geq 2$, for each $i \in I$, let $b_{i,k} = a_{i,k} \equiv a_i$ (*i.e.* constant symmetric inertial parameters), then (2.3) tells us that the $a_i$ must live in a ball,

$$\Big(\frac{1}{2\underline{\gamma}\mu} + \frac{1}{2\nu/L^2}\Big) \sum_{i\in I} a_i^2 < \underline{\beta}.$$

**An empirical approach for inertial parameters** Besides Theorem 2.2, we also provide an empirical bound for the choice of the inertial parameters. Consider the setting: $\gamma_k \equiv \gamma \in ]0, 1/L[$ and $b_{i,k} = a_{i,k} \equiv a_i \in ]-1, 2[, i \in I$. We have the following empirical bound for the summand $\sum_{i\in I} a_i$:

$$\textstyle\sum_i a_i \in \big]0, \min\big\{1, \frac{1/L-\gamma}{|2\gamma-1/L|}\big\}\big[. \tag{2.4}$$

To ensure the convergence $\{x_k\}_{k\in\mathbb{N}}$, an online updating rule should be applied together with the empirical bound. More precisely, choose $a_i$ according to (2.4). Then for each $k \in \mathbb{N}$, let $b_{i,k} = a_{i,k}$ and choose $a_{i,k}$ such that $\sum_i a_{i,k} = \min\{\sum_i a_i, c_k\}$ where $c_k > 0$ is such that $\{c_k \sum_{i\in I} \|x_{k-i} - x_{k-i-1}\|\}_{k\in\mathbb{N}}$ is summable. For instance,

$$c_k = \frac{c}{k^{1+q}\sum_{i\in I}\|x_{k-i} - x_{k-i-1}\|}, \ c > 0, \ q > 0.$$

**Remark 2.4.** The allowed choices of the summand $\sum_i a_i$ by (2.4) is larger than those of Theorem 2.2. For instance, (2.4) allows $\sum_i a_i = 1$ for $\gamma \in ]0, \frac{2}{3L}]$. While for Theorem 2.2, $\sum_i a_i = 1$ can be reached only when $\gamma \to 0$.

# 3 Local convergence properties of MiFB

Here, we present the local convergence analysis of the proposed method under partial smoothness.

## 3.1 Partial smoothness

Let $\mathcal{M} \subset \mathbb{R}^n$ be a $C^2$-smooth submanifold, let $\mathcal{T}_{\mathcal{M}}(x)$ the tangent space of $\mathcal{M}$ at any point $x \in \mathcal{M}$.

**Definition 3.1.** The function $R : \mathbb{R}^n \to \mathbb{R} \cup \{+\infty\}$ is $C^2$*-partly smooth at* $\bar{x} \in \mathcal{M}$ *relative to* $\mathcal{M}$ *for* $\bar{v} \in \partial R(\bar{x}) \neq \emptyset$ if $\mathcal{M}$ is a $C^2$-submanifold around $\bar{x}$, and

(i) (Smoothness): $R$ restricted to $\mathcal{M}$ is $C^2$ around $\bar{x}$;

(ii) (Regularity): $R$ is regular at all $x \in \mathcal{M}$ near $\bar{x}$ and $R$ is $r$-prox-regular at $\bar{x}$ for $\bar{v}$;

(iii) (Sharpness): $\mathcal{T}_{\mathcal{M}}(\bar{x}) = \mathrm{par}(\partial R(x))^{\perp}$;

(iv) (Continuity): The set-valued mapping $\partial R$ is continuous at $\bar{x}$ relative to $\mathcal{M}$.

We denote the class of partly smooth functions at $x$ relative to $\mathcal{M}$ for $v$ as $\mathrm{PSF}_{x,v}(\mathcal{M})$. Partial smoothness was first introduced in [18] and its directional version stated here is due to [21, 15]. Prox-regularity is sufficient to ensure that the partly smooth submanifolds are locally unique [21, Corollary 4.12], [15, Lemma 2.3 and Proposition 10.12].

### 3.2 Finite activity identification

One of the key consequences of partial smoothness is finite identification of the partial smoothness submanifold associated to $R$ for problem ($\mathcal{P}$). This is formalized in the following statement.

**Theorem 3.2 (Finite activity identification).** *Suppose that Algorithm 1 is run under the conditions of Theorem 2.2, such that the generated sequence $\{x_k\}_{k\in\mathbb{N}}$ converges to a critical point $x^\star \in \mathrm{crit}(\Phi)$. Assume that $R \in \mathrm{PSF}_{x^\star,-\nabla F(x^\star)}(\mathcal{M}_{x^\star})$ and the non-degeneracy condition*

$$-\nabla F(x^\star) \in \mathrm{ri}\big(\partial R(x^\star)\big), \tag{ND}$$

*holds. Then, $x_k \in \mathcal{M}_{x^\star}$ for all $k$ large enough.*

See the Section B for the proof. This result generalises that of [23] to the non-convex case and multiple inertial steps.

### 3.3 Local linear convergence

Given $\gamma \in ]0, \frac{1}{L}[$ and a critical point $x^\star \in \mathrm{crit}(\Phi)$, let $\mathcal{M}_{x^\star}$ be a $C^2$-smooth submanifold and $R \in \mathrm{PSF}_{x^\star,-\nabla F(x^\star)}(\mathcal{M}_{x^\star})$. Denote $T_{x^\star} \stackrel{\text{def}}{=} \mathcal{T}_{\mathcal{M}_{x^\star}}(x^\star)$ and the following matrices which are all symmetric,

$$H \stackrel{\text{def}}{=} \gamma \mathrm{P}_{T_{x^\star}} \nabla^2 F(x^\star) \mathrm{P}_{T_{x^\star}}, \quad G \stackrel{\text{def}}{=} \mathrm{Id} - H, \quad Q \stackrel{\text{def}}{=} \gamma \nabla^2_{\mathcal{M}_{x^\star}} \Phi(x^\star) \mathrm{P}_{T_{x^\star}} - H, \tag{3.1}$$

where $\nabla^2_{\mathcal{M}_{x^\star}} \Phi$ is the Riemannian Hessian of $\Phi$ along the submanifold $\mathcal{M}_{x^\star}$ (readers may refer to the supplementary material from more details on differential calculus on Riemannian manifolds).

To state our local linear convergence result, the following assumptions will play a key role.

**Restricted injectivity** Besides the local $C^2$-smoothness assumption on $F$, following the idea of [22, 23], we assume the restricted injectivity condition,

$$\ker\big(\nabla^2 F(x^\star)\big) \cap T_{x^\star} = \{0\}. \tag{RI}$$

**Positive semi-definiteness of $Q$** Assume that $Q$ is *positive semi-definite*, *i.e.* $\forall h \in T_{x^\star}$,

$$\langle h,\, Qh\rangle \geq 0. \tag{3.2}$$

Under (3.2), $\mathrm{Id} + Q$ is symmetric positive definite, hence invertible, we denote $P \stackrel{\text{def}}{=} (\mathrm{Id} + Q)^{-1}$.

**Convergent parameters** The parameters of MiFB (Algorithm 1), are convergent, *i.e.*

$$a_{i,k} \to a_i,\; b_{i,k} \to b_i,\, \forall i \in I \;\text{ and }\; \gamma_k \to \gamma \in [\underline{\gamma}, \min\{\bar{\gamma}, \bar{r}\}], \tag{3.3}$$

where $\bar{r} < r$, and $r$ is the prox-regularity modulus of $R$ (see Definition 3.1).

**Remark 3.3.**

(i) Condition (3.2) can be met by various non-convex functions, such as polyhedral functions, including the $\ell_0$ pseudo-norm discussed in Example 1.1 and Section 4.1.

(ii) Condition (3.3) asserts that both the inertial parameters $(a_{i,k}, b_{i,k})$ and the step-size $\gamma_k$ should converge to some limit points, and this condition cannot be relaxed in general.

(iii) It can be shown that conditions (3.2) and (RI) together imply that $x^\star$ is a local minimizer of $\Phi$ in ($\mathcal{P}$), and $\Phi$ grows at least quadratically near $x^\star$. The arguments to prove this are essentially adapted from those used to show [23, Proposition 4.1(ii)].

We need the following notations:

$$
\begin{gathered}
M_0 \stackrel{\text{def}}{=} (a_0 - b_0)P + (1+b_0)PG,\ \ M_s \stackrel{\text{def}}{=} -(a_{s-1} - b_{s-1})P - b_{s-1}PG, \\
M_i \stackrel{\text{def}}{=} -\big((a_{i-1} - a_i) - (b_{i-1} - b_i)\big)P - (b_{i-1} - b_i)PG,\ \ i = 1, ..., s-1, \\
M \stackrel{\text{def}}{=} \begin{bmatrix} M_0 & \cdots & M_{s-1} & M_s \\ \mathrm{Id} & \cdots & 0 & 0 \\ \vdots & \ddots & \vdots & \vdots \\ 0 & \cdots & \mathrm{Id} & 0 \end{bmatrix},\ \ d_k \stackrel{\text{def}}{=} \begin{pmatrix} x_k - x^\star \\ \vdots \\ x_{k-s} - x^\star \end{pmatrix}.
\end{gathered} \tag{3.4}
$$

**Theorem 3.4 (Local linear convergence).** *Suppose that the MiFB Algorithm 1 is run under the setting of Theorem 3.2. Moreover, assume that (RI), (3.2) and (3.3) hold. Then for all $k$ large enough,*

$$d_{k+1} = M d_k + o(\|d_k\|). \tag{3.5}$$

*If $\rho(M) < 1$, then given any $\rho \in ]\rho(M), 1[$, there exists $K \in \mathbb{N}$ such that $\forall k \geq K$,*

$$\|x_k - x^\star\| = O(\rho^{k-K}). \tag{3.6}$$

*In particular, if $s = 1$, then $\rho(M) < 1$ if $R$ is locally polyhedral around $x^\star$ or if $a_0 = b_0$.*

See the Section C for the proof.

**Remark 3.5.**

(i) When $s = 1$, $\rho(M)$ can be given explicitly in terms of the parameters of the algorithm (*i.e.* $a_0$, $b_0$ and $\gamma$), see [23, Section 4.2] for details. However, the spectral analysis of $M$ becomes much more complicated to get for $s \geq 2$, where the analysis of at least cubic equations are involved. Therefore, for the sake of brevity, we shall skip the detailed discussion here.

(ii) When $s = 1$, it was shown in [23] that the optimal convergence rate that can be obtained by 1-step inertial scheme with fixed $\gamma$ is $\rho^\star_{s=1} = 1 - \sqrt{1 - \tau\gamma}$, where from condition (RI), continuity of $\nabla^2 F$ at $x^\star$ implies that there exists $\tau > 0$ and a neighbourhood of $x^\star$ such that $\langle h,\, \nabla^2 F(x^\star)h\rangle \geq \tau\|h\|^2$, for all $h \in T_{x^\star}$. As we will see in the numerical experiments of Section 4, such a rate can be improved by our multi-step inertial scheme. Taking $s = 2$ for example, we will show that for a certain class of functions, the optimal local linear rate is close to or even is $\rho^\star_{s=2} = 1 - \sqrt[3]{1 - \tau\gamma}$, which is obviously faster than $\rho^\star_{s=1}$.

(iii) Though it can be satisfied for many problems in practice, the restricted injectivity (RI) can be removed for some penalties $R$, for instance, when $R$ is locally polyhedral near $x^\star$.

# 4 Numerical experiments

In this section, we illustrate our results with some numerical experiments carried out on the problems in Example 1.1, 1.2 and 1.3.

### 4.1 Examples of KL and partly smooth functions

All the objectives $\Phi$ in the above mentioned examples are continuous KL functions. Indeed, in Example 1.1 and 1.2, $\Phi$ is the sum of semi-algebraic functions which is also semi-algebraic. In Example 1.3, $\Phi$ is also algebraic when $G$ is the squared hinge loss, and definable in an o-minimal structure for the logistic loss (see *e.g.* [31] for material on o-minimal structures).

Moreover, $R$ is partly smooth in all these examples as we show now.

**Example 4.1 ($\ell_0$ pseudo-norm).** The $\ell_0$ pseudo-norm is locally constant. Moreover, it is regular on $\mathbb{R}^n$ ([16, Remark 2]) and its subdifferential is given by (see [16, Theorem 1])

$$\partial\|x\|_0 = \mathrm{span}\big((e_i)_{i\in\mathrm{supp}(x)^c}\big),$$

where $(e_i)_{i=1,\ldots,n}$ is the standard basis, and $\mathrm{supp}(x) = \big\{i: \ x_i \neq 0\big\}$. The proximity operator of $\ell_0$-norm is given by hard-thresholding,

$$\mathrm{prox}_{\gamma\|x\|_0}(z) = \begin{cases} z & \text{if } |z| > \sqrt{2\gamma}, \\ \mathrm{sign}(z)[0,z] & \text{if } |z| = \sqrt{2\gamma}, \\ 0 & \text{if } |z| < \sqrt{2\gamma}. \end{cases}$$

It can then be easily verified that the $\ell_0$ pseudo-norm is partly smooth at any $x$ relative to the subspace

$$\mathcal{M}_x = T_x = \big\{z \in \mathbb{R}^n: \ \mathrm{supp}(z) \subset \mathrm{supp}(x)\big\}.$$

It is also prox-regular at $x$ for any bounded $v \in \partial\|x\|_0$. Note also condition (ND) is automatically verified and that the Riemannian gradient and Hessian along $T_x$ of $\|\cdot\|_0$ vanish.

**Example 4.2 (Rank).** The rank function is the spectral extension of $\ell_0$ pseudo-norm to matrix-valued data $x \in \mathbb{R}^{n_1\times n_2}$ [20]. Consider a singular value decomposition (SVD) of $x$, *i.e.* $x = U\mathrm{diag}(\sigma(x))V^*$, where $U = \{u_1,\ldots,u_n\}, V = \{v_1,\ldots,v_n\}$ are orthonormal matrices, and $\sigma(x) = (\sigma_i(x))_{i=1,\ldots,n}$ is the vector of singular values. By definition, $\mathrm{rank}(x) \stackrel{\text{def}}{=} \|\sigma(x)\|_0$. Thus the rank function is partly smooth relative at $x$ to the set of fixed rank matrices

$$\mathcal{M}_x = \big\{z \in \mathbb{R}^{n_1\times n_2}: \ \mathrm{rank}(z) = \mathrm{rank}(x)\big\},$$

which is a $C^2$-smooth submanifold [19]. The tangent space of $\mathcal{M}_x$ at $x$ is

$$\mathcal{T}_{\mathcal{M}_x}(x) = T_x = \big\{z \in \mathbb{R}^{n_1\times n_2}: \ u_i^* z v_j = 0, \text{for all } r < i \leq n_1, \ r < j \leq n_2\big\},$$

The $\mathrm{rank}$ function is also regular its subdifferential reads

$$\partial\mathrm{rank}(x) = U\partial\big(\|\sigma(x)\|_0\big)V^* = U\mathrm{span}\big((e_i)_{i\in\mathrm{supp}(\sigma(x))^c}\big)V^*,$$

which is a vector space (see [16, Theorem 4 and Proposition 1]). The proximity operator of $\mathrm{rank}$ function amounts to applying hard-thresholding to the singular values. Observe that by definition of $\mathcal{M}_x$, the Riemannian gradient and Hessian of the rank function along $\mathcal{M}_x$ also vanish.

For Example 1.2, it is worth noting from the above examples and separability of the regularizer that the latter is also partly smooth relative to the cartesian product of the partial smoothness submanifolds of $\ell_0$ and the rank function.

## 4.2 Experimental results

For the problem in Example 1.1, we generated $y = Ax_{\text{ob}} + \omega$ with $m = 48$, $n = 128$, the entries of $A$ are i.i.d. zero-mean and unit variance Gaussian, $x_{\text{ob}}$ is 8-sparse, and $\omega \in \mathbb{R}^m$ is an additive noise with small variance.

For the problem in Example 1.2, we generated $y = x_s + x_l + \omega$, with $n_1 = n_2 = 50$, $x_s$ is 250-sparse, and the rank of $x_l$ is 5, and $\omega$ is an additive noise with small variance.

For Example 1.3, we generated $m = 64$ training samples with $n = 96$-dimensional feature space.

For all presented numerical results, 3 different settings were tested:

- the FB method, with $\gamma_k \equiv 0.3/L$, noted as "FB";
- MiFB with $s = 1$, $b_k = a_k \equiv a$ and $\gamma_k \equiv 0.3/L$, noted as "1-iFB";
- MiFB with $s = 2$, $b_{i,k} = a_{i,k} \equiv a_i, i = 0, 1$ and $\gamma_k \equiv 0.3/L$, noted as "2-iFB".

**Tightness of theoretical prediction** The convergence profiles of $\|x_k - x^\star\|$ are shown in Figure 1. As it can be seen from all the plots, finite identification and local linear convergence indeed occur. The positions of the *green dots* indicate the iteration from which $x_k$ numerically identifies the submanifold $\mathcal{M}_{x^\star}$. The solid lines ("P") represents practical observations, while the dashed lines ("T") denotes theoretical predictions.

As the Riemannian Hessians of $\ell_0$ and the rank both vanish in all examples, our predicted rates coincide exactly with the observed ones (same slopes for the dashed and solid lines).

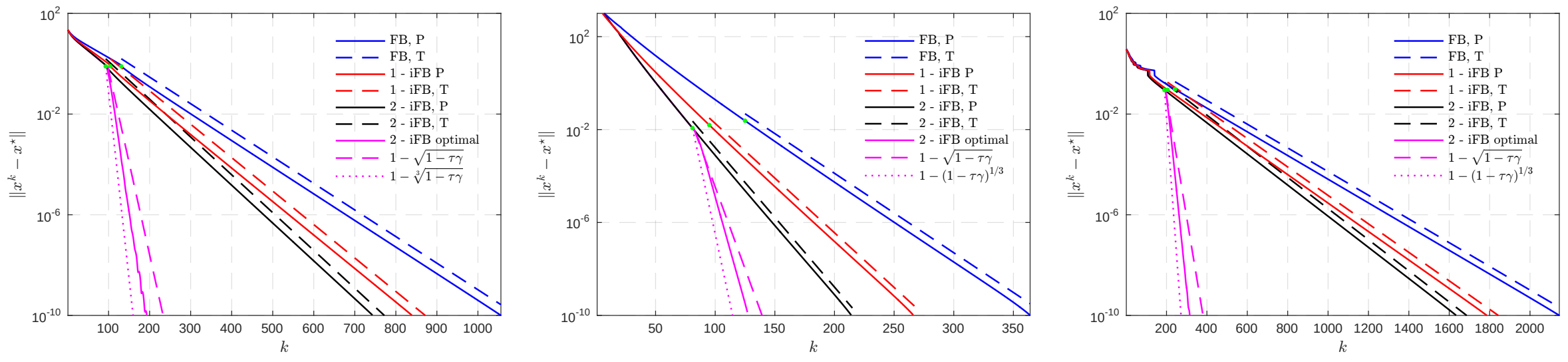

(a) Sparse regression (b) PCP (c) Sparse SVM

Figure 1: Finite identification and local linear convergence of MiFB under different inertial settings in terms of $\|x_k - x^\star\|$. "P" stands for practical observation and "T" indicates the theoretical estimate. We fix $\gamma_k \equiv 0.3/L$ for all tests. For the 2 inertial schemes, inertial parameters are first chosen such that (2.3) holds. The position of the green dot in each plot indicates the iteration beyond which identification of $\mathcal{M}_{x^\star}$ occurs.

**Comparison of the methods** Under the tested settings, we draw the following remarks on the comparison of the inertial schemes:

- The MiFB scheme is much faster than FB both globally and locally. Finite activity identification also occurs earlier for MiFB than for FB;
- Comparing the two MIFB inertial schemes, "2-iFB" outperforms "1-iFB", showing the advantages of a 2-step inertial scheme over the 1-step one.

**Optimal first-order method** To highlight the potential of multiple steps in MiFB, for the "2-iFB" scheme, we also added an example where we locally optimized the rate for the inertial parmeters. See the *magenta* lines all the examples, where the solid line corresponds to the observed profile for the optimal inertial parameters, the *dashed* line stands for the rate $1 - \sqrt{1 - \tau\gamma}$, and the *dotted* line is that of $1 - \sqrt[3]{1 - \tau\gamma}$, which shows indeed that a faster linear rate can be obtained owing to multiple inertial parameters.

We refer to [23, Section 4.5] for the optimal choice of inertial parameters for the case $s = 1$.

### 4.3 The empirical bound (2.4) and inertial steps $s$

We now present a short comparison of the empirical bound (2.4) of inertial parameters and different choices of $s$ under bigger choice of $\gamma = 0.8/L$. MiFB with 3 inertial steps, *i.e.* $s = 3$, is added which is noted as "3-iFB", see the *magenta* line in Figure 2.

Similar to the above experiments, we choose $b_{i,k} = a_{i,k} \equiv a_i, i \in I$, and "Thm 2.2" means that $a_i$'s are chosen according to Theorem 2.2, while "Bnd (2.4)" means that $a_i$'s are chosen based on the empirical bound (2.4). We can infer from Figure 2 the following.

- Compared to the results in Figure 1, a bigger choice of $\gamma$ leads to faster convergence. Yet still, under the same choice of $\gamma$, MiFB is faster than FB both locally and globally;
- For either "Thm 2.2" or "Bnd (2.4)", the performance of the three MiFB schemes are close, this is mainly due to the fact that values of the sum $\sum_{i\in I} a_i$ for each scheme are close.
- Then between "Thm 2.2" and "Bnd (2.4)", "Bnd (2.4)" shows faster convergence result, since the allowed value of $\sum_{i\in I} a_i$ of (2.4) is bigger than that of Theorem 2.2.

It should be noted that, when $\gamma \in ]0, \frac{2}{3L}]$, the largest value of $\sum_{i\in I} a_i$ allowed by (2.4) is $1$. If we choose $\sum_{i\in I} a_i$ equal or very close to $1$, then it can be observed in practice that MiFB locally oscillates, which is a well-known property of the FISTA scheme [6, 11]. We refer to [23, Section 4.4] for discussions of the properties of such oscillation behaviour.

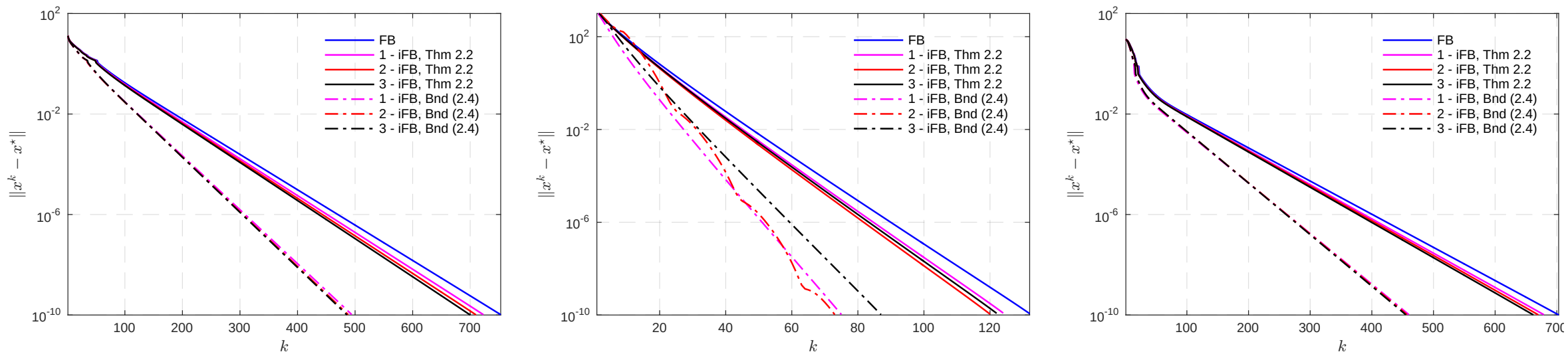

(a) Sparse regression (b) PCP (c) Sparse SVM

Figure 2: Comparison of MiFB under different inertial settings. We fix $\gamma_k \equiv 0.8/L$ for all tests. For the three inertial schemes, the inertial parameters were chosen such that (2.3) holds.

## A Proof of Theorem 2.2

**Lemma A.1.** *Let $\{d_k\}_{k\in\mathbb{N}}, \{\delta_k\}_{k\in\mathbb{N}}$ be two non-negative sequences, and $\omega \in \mathbb{R}^s$ such that*

$$d_{k+1} \leq \sum\nolimits_{i\in I} \omega_i d_{k-i} + \delta_k, \tag{A.1}$$

*for all $k \geq s$. If $\sum_i \omega_i \in [0, 1[$ and $\sum_{k\in\mathbb{N}} \delta_k < +\infty$, then*

$$\sum\nolimits_{k\in\mathbb{N}} d_k < +\infty.$$

**Remark A.2.** Lemma A.1 is an extension of [9, Lemma 3]. It should be noted that in our case, non-negativity is *not* imposed to the weight $\omega_i$'s, but only the sum of them. In fact, we can even afford all $\omega_i$'s to be negative, as long as $\sum_{i\in I} \omega_i d_{k-i} + \delta_k$ is positive for all $k \in \mathbb{N}$.

**Proof.** From (A.1), suppose that $d_{-1} = d_{-2} = d_{-s+1} = 0$, then sum up for both sides from $k = 0$,

$$\begin{aligned}\sum_{k\in\mathbb{N}} d_{k+1} \le \sum_{k\in\mathbb{N}}\sum_{i\in I}\omega_i d_{k-i} + \sum_{k\in\mathbb{N}}\delta_k &\Longrightarrow \sum_{k\in\mathbb{N}} d_k \le d_0 + \sum_{i\in I}\omega_i \sum_{k\in\mathbb{N}} d_k + \sum_{k\in\mathbb{N}}\delta_k \\ &\Longrightarrow \Big(1 - \sum_{i\in I}\omega_i\Big)\sum_{k\in\mathbb{N}} d_k \le d_0 + \sum_{k\in\mathbb{N}}\delta_k.\end{aligned}$$

Since we assume $\sum_{i\in I}\omega_i < 1$ and $\delta_k$ is summable, then we have

$$\sum_{k\in\mathbb{N}} d_k \le \Big(1 - \sum_{i\in I}\omega_i\Big)^{-1}\Big(d_0 + \sum_{k\in\mathbb{N}}\delta_k\Big) < +\infty,$$

which concludes the proof. □

Define $\Delta_k \stackrel{\text{def}}{=} \|x_k - x_{k-1}\|$.

**Lemma A.3.** *For the update of $x_{k+1}$ in (1.7), given any $k \in \mathbb{N}$, define*

$$g_{k+1} \stackrel{\text{def}}{=} \frac{1}{\gamma_k}(y_{a,k} - x_{k+1}) - \nabla F(y_{b,k}) + \nabla F(x_{k+1}).$$

*We have $g_{k+1} \in \partial\Phi(x_{k+1})$, and moreover,*

$$\|g_{k+1}\| \le \Big(\frac{1}{\underline{\gamma}} + L\Big)\Delta_{k+1} + \sum_{i\in I}\Big(\frac{|a_{i,k}|}{\underline{\gamma}} + |b_{i,k}|\Big)\Delta_{k-i}. \tag{A.2}$$

**Proof.** From the definition of proximity operator and the update of $x_{k+1}$ (1.7), we have $y_{a,k} - \gamma_k\nabla F(y_{b,k}) - x_{k+1} \in \gamma_k\partial R(x_{k+1})$, add $\gamma_k\nabla F(x_{k+1})$ to both sides, then

$$g_{k+1} = \frac{y_{a,k} - \gamma_k\nabla F(y_{b,k}) - x_{k+1} + \gamma_k\nabla F(x_{k+1})}{\gamma_k} \in \partial\Phi(x_{k+1}).$$

Then, apply the triangle inequality and the Lipschitz continuity of $\nabla F$, we get

$$\begin{aligned}\|g_{k+1}\| &= \|\tfrac{1}{\gamma_k}(y_{a,k} - x_{k+1}) - \nabla F(y_{b,k}) + \nabla F(x_{k+1})\| \\ &\le \frac{1}{\gamma_k}\|y_{a,k} - x_{k+1}\| + L\|y_{b,k} - x_{k+1}\| \\ &\le \frac{1}{\gamma_k}\Big(\Delta_{k+1} + \sum_{i\in I}|a_{i,k}|\Delta_{k-i}\Big) + L\Big(\Delta_{k+1} + \sum_{i\in I}|b_{i,k}|\Delta_{k-i}\Big) \\ &\le \Big(\frac{1}{\underline{\gamma}} + L\Big)\Delta_{k+1} + \sum_{i\in I}\Big(\frac{|a_{i,k}|}{\underline{\gamma}} + |b_{i,k}|\Big)\Delta_{k-i},\end{aligned}$$

which concludes the proof. □

**Lemma A.4.** *For Algorithm 1, given the parameters $\gamma_k, a_{i,k}, b_{i,k}$, the following inequality holds*

$$\Phi(x_{k+1}) + \underline{\beta}\Delta_{k+1}^2 \le \Phi(x_k) + \sum_{i\in I}\bar{\alpha}_i\Delta_{k-i}^2. \tag{A.3}$$

**Proof.** Define the function

$$\mathcal{L}_k(x) = \gamma_k R(x) + \frac{1}{2}\|x - y_{a,k}\|^2 + \gamma_k\langle x,\, \nabla F(y_{b,k})\rangle.$$

It can be shown that the update of $x_{k+1}$ in (1.7) is equivalent to

$$x_{k+1} \in \operatorname{argmin}_{x\in\mathbb{R}^n}\mathcal{L}_k(x), \tag{A.4}$$

which means that $\mathcal{L}_k(x_{k+1}) \leq \mathcal{L}_k(x_k)$, and

$$R(x_{k+1}) + \frac{1}{2\gamma_k}\|x_{k+1} - y_{a,k}\|^2 + \langle x_{k+1},\, \nabla F(y_{b,k})\rangle \leq R(x_k) + \frac{1}{2\gamma_k}\|x_k - y_{a,k}\|^2 + \langle x_k,\, \nabla F(y_{b,k})\rangle,$$

and leads to,

$$\begin{aligned} R(x_k) &\geq R(x_{k+1}) + \frac{1}{2\gamma_k}\|x_{k+1} - y_{a,k}\|^2 + \langle x_{k+1} - x_k,\, \nabla F(y_{b,k})\rangle - \frac{1}{2\gamma_k}\|x_k - y_{a,k}\|^2 \\ &= R(x_{k+1}) + \langle x_{k+1} - x_k,\, \nabla F(x_k)\rangle + \frac{1}{2\gamma_k}\Delta_{k+1}^2 \\ &\quad + \frac{1}{\gamma_k}\langle x_k - x_{k+1},\, \textstyle\sum_{i\in I} a_{i,k}(x_{k-i} - x_{k-i-1})\rangle + \langle x_{k+1} - x_k,\, \nabla F(y_{b,k}) - \nabla F(x_k)\rangle. \end{aligned} \tag{A.5}$$

Since $F$ is $L$-Lipschitz, then

$$\langle \nabla F(x_k),\, x_{k+1} - x_k\rangle \geq F(x_{k+1}) - F(x_k) - \frac{L}{2}\Delta_{k+1}^2.$$

Apply Young's inequality, we obtain

$$\begin{aligned} \langle x_k - x_{k+1},\, \textstyle\sum_{i\in I} a_{i,k}(x_{k-i} - x_{k-i-1})\rangle &\geq -\Big(\frac{\mu}{2}\Delta_{k+1}^2 + \frac{1}{2\mu}\|\textstyle\sum_{i\in I} a_{i,k}(x_{k-i} - x_{k-i-1})\|^2\Big) \\ &\geq -\Big(\frac{\mu}{2}\Delta_{k+1}^2 + \textstyle\sum_{i\in I}\frac{s a_{i,k}^2}{2\mu}\Delta_{k-i}^2\Big), \end{aligned} \tag{A.6}$$

where $\mu > 0$. Then similarly, for $\nu > 0$, we have

$$\begin{aligned} \langle x_{k+1} - x_k,\, \nabla F(y_{b,k}) - \nabla F(x_k)\rangle &\geq -\Big(\frac{\nu}{2}\Delta_{k+1}^2 + \frac{1}{2\nu}\|\nabla F(y_{b,k}) - \nabla F(x_k)\|^2\Big) \\ &\geq -\Big(\frac{\nu}{2}\Delta_{k+1}^2 + \textstyle\sum_{i\in I}\frac{s b_{i,k}^2 L^2}{2\nu}\Delta_{k-i}^2\Big). \end{aligned} \tag{A.7}$$

Combining the above 3 inequalities with (A.5) leads to

$$\Phi(x_{k+1}) + \beta_k\Delta_{k+1}^2 \leq \Phi(x_k) + \textstyle\sum_{i\in I}\Big(\frac{s a_{i,k}^2}{2\gamma_k\mu} + \frac{s b_{i,k}^2 L^2}{2\nu}\Big)\Delta_{k-i}^2 = \Phi(x_k) + \sum_{i\in I}\alpha_{k,i}\Delta_{k-i}^2. \tag{A.8}$$

Therefore, we obtain

$$\Phi(x_{k+1}) + \underline{\beta}\Delta_{k+1}^2 \leq \Phi(x_{k+1}) + \beta_k\Delta_{k+1}^2 \leq \Phi(x_k) + \textstyle\sum_{i\in I}\alpha_{k,i}\Delta_{k-i}^2 \leq \Phi(x_k) + \sum_{i\in I}\bar{\alpha}_i\Delta_{k-i}^2,$$

which concludes the proof. □

Define $\mathbb{R}_s^n$ the product space $\mathbb{R}_s^n \stackrel{\text{def}}{=} \underbrace{\mathbb{R}^n \times \cdots \times \mathbb{R}^n}_{s \text{ times}}$ and $z_k = (x_k, x_{k-1}, ..., x_{k-s+1}) \in \mathbb{R}_s^n$. Then given $z_k$, define the function

$$\Psi(z_k) = \Phi(x_k) + \textstyle\sum_{i\in I}\sum_{j=i}^{s-1}\bar{\alpha}_j\Delta_{k-i}^2,$$

which is is a KL function if $\Phi$ is. Denote $\mathcal{C}_{x_k}, \mathcal{C}_{z_k}$ the set of cluster points of sequences $x_k$ and $z_k$ respectively, and $\text{crit}(\Psi) = \{z = (x, ..., x) \in \mathbb{R}_s^n : x \in \text{crit}(\Phi)\}$.

**Lemma A.5.** *For Algorithm 1, choose $\mu, \nu, \gamma_k, a_{i,k}, b_{i,k}$ such that (2.3) holds. If $\Phi$ is bounded from below, then*

(i) $\sum_{k\in\mathbb{N}} \Delta_k^2 < +\infty$;

(ii) *The sequence* $\Psi(z_k)$ *is monotonically decreasing and convergent;*

(iii) *The sequence* $\Phi(x_k)$ *is convergent.*

**Proof.** Define

$$\delta = \underline{\beta} - \sum_{i\in I}\bar{\alpha}_i > 0.$$

From the Lemma A.4, we have

$$\delta\Delta_{k+1}^2 \le \big(\Phi(x_k) - \Phi(x_{k+1})\big) + \sum_{i\in I}\bar{\alpha}_i(\Delta_{k-i}^2 - \Delta_{k+1}^2).$$

Since we let $x_{1-s} = \ldots = x_0 = x_1$, for the above inequality, sum over $k$ we get

$$\begin{aligned}\delta \sum_{k\in\mathbb{N}} \Delta_{k+1}^2 &\le \sum_{k\in\mathbb{N}}\big(\Phi(x_k) - \Phi(x_{k+1})\big) + \sum_{k\in\mathbb{N}}\sum_{i\in I}\bar{\alpha}_i(\Delta_{k-i}^2 - \Delta_{k+1}^2)\\ &\le \Phi(x_0) + \sum_{i\in I}\bar{\alpha}_i\sum_{k\in\mathbb{N}}(\Delta_{k-i}^2 - \Delta_{k+1}^2)\\ &= \Phi(x_0) + \sum_{i\in I}\bar{\alpha}_i\sum_{j=1-i}^{1}\Delta_j^2 = \Phi(x_0),\end{aligned}$$

which means, as $\Phi(x_0)$ is bounded,

$$\sum_{k\in\mathbb{N}}\Delta_{k+1}^2 \le \frac{\Phi(x_0)}{\delta} < +\infty.$$

From Lemma A.4, by pairing terms on both sides of (A.3), we get

$$\Psi(z_{k+1}) + \big(\underline{\beta} - \sum_{i\in I}\bar{\alpha}_i\big)\Delta_{k+1}^2 \le \Psi(z_k).$$

Since we assume $\underline{\beta} - \sum_{i\in I}\bar{\alpha}_i > 0$, hence $\Psi(z_k)$ is monotonically non-increasing. The convergence of $\Phi(x_k)$ is straightforward. □

**Lemma A.6.** *For Algorithm 1, choose* $\mu, \nu, \gamma_k, a_{i,k}, b_{i,k}$ *such that* (2.3) *holds. If* $\Phi$ *is bounded from below and* $\{x_k\}_{k\in\mathbb{N}}$ *is bounded, then* $x_k$ *converges to a critical point of* $\Phi$.

**Proof.** Since $\{x_k\}_{k\in\mathbb{N}}$ is bounded, there exists a subsequence $\{x_{k_j}\}_{k\in\mathbb{N}}$ and cluster point $\bar{x}$ such that $x_{k_j} \to \bar{x}$ as $j \to \infty$. Next we show that $\Phi(x_{k_j}) \to \Phi(\bar{x})$ and that $\bar{x}$ is a critical point of $\Phi$.

Since $R$ is lsc, then $\liminf_{j\to\infty} R(x_{k_j}) \ge R(\bar{x})$. From (A.4), we have $\mathcal{L}_{k_j-1}(x_{k_j}) \le \mathcal{L}_{k_j-1}(\bar{x})$,

$$\begin{aligned}R(\bar{x}) &\ge R(x_{k_j}) + \frac{1}{2\gamma_{k_j-1}}\|x_{k_j} - y_{a,k_j-1}\|^2 + \langle x_{k_j} - \bar{x},\, \nabla F(y_{b,k_j-1})\rangle - \frac{1}{2\gamma_{k_j-1}}\|\bar{x} - y_{a,k_j-1}\|^2\\ &= R(x_{k_j}) + \frac{1}{2\gamma_{k_j-1}}(\|x_{k_j} - \bar{x}\|^2 + 2\langle x_{k_j} - \bar{x},\, \bar{x} - y_{a,k_j-1}\rangle) + \langle x_{k_j} - \bar{x},\, \nabla F(y_{b,k_j-1})\rangle\end{aligned}$$

Since $\Delta_k^2 \to 0$ and $x_{k_j} \to \bar{x}$, then taking the above inequality to limit we have $\limsup_{j\to\infty} R(x_{k_j}) \le R(\bar{x})$. As a result, $\lim_{k\to\infty} R(x_{k_j}) = R(\bar{x})$. Since $F$ is continuous, then $F(x_{k_j}) \to F(\bar{x})$, hence $\Phi(x_{k_j}) \to \Phi(\bar{x})$.

Furthermore, owing to Lemma A.3, $g_{k_j} \in \partial\Phi(x_{k_j})$, and (i) of Lemma A.5 we have $g_{k_j} \to 0$ as $k \to \infty$. As a consequence,

$$g_{k_j} \in \partial\Phi(x_{k_j}),\ \ (x_{k_j}, g_{k_j}) \to (\bar{x}, 0)\ \text{ and }\ \Phi(x_{k_j}) \to \Phi(\bar{x}),$$

as $j \to \infty$. Hence $0 \in \partial\Phi(\bar{x})$, *i.e.* $\bar{x}$ is a critical point. □

Now we present the proof of Theorem 2.2.

**Proof of Theorem 2.2.** Putting together the above lemmas, we draw the following useful conclusions:

**(R.1)** Denote $\delta = \underline{\beta} - \sum_{i\in I} \bar{\alpha}_i$, then $\Psi(z_{k+1}) + \delta\Delta_{k+1}^2 \le \Psi(z_k)$;

**(R.2)** Define

$$w_{k+1} \overset{\text{def}}{=} \begin{pmatrix} g_{k+1} + 2\sum_{i=0}^{s-1} \bar{\alpha}_i(x_{k+1} - x_k) \\ 2\sum_{i=0}^{s-1} \bar{\alpha}_i(x_k - x_{k+1}) + 2\sum_{i=1}^{s-1} \bar{\alpha}_i(x_k - x_{k-1}) \\ \vdots \\ 2\bar{\alpha}_{s-1}(x_{k+2-s} - x_{k+1-s}) \end{pmatrix},$$

then we have $w_{k+1} \in \partial\Psi(z_{k+1})$. Owing to Lemma A.3, there exists a $\sigma > 0$ such that $\|w_{k+1}\| \le \sigma\sum_{j=k+2-s}^{k+1} \Delta_j$;

**(R.3)** if $x_{k_j}$ is a subsequence such that $x_{k_j} \to \bar{x}$, then $\Psi(z_k) \to \Psi(\bar{z})$ where $\bar{z} = (\bar{x}, ..., \bar{x})$.

**(R.4)** $\mathcal{C}_{z_k} \subseteq \mathrm{crit}(\Psi)$;

**(R.5)** $\lim_{k\to\infty} \mathrm{dist}(z_k, \mathcal{C}_{z_k}) = 0$;

**(R.6)** $\mathcal{C}_{z_k}$ is non-empty, compact and connected;

**(R.7)** $\Psi$ is finite and constant on $\mathcal{C}_{z_k}$.

Next we prove the claims of Theorem 2.2.

(i) Consider a critical point of $\Phi$, $\bar{x} \in \mathrm{crit}(\Phi)$, such that $\bar{z} = (\bar{x}, ..., \bar{x}) \in \mathcal{C}_{z_k}$, then owing to (**R.3**), we have $\Psi(z_k) \to \Psi(\bar{z})$.

Suppose there exists $K$ such that $\Psi(z_K) = \Psi(\bar{z})$, then the descent property (**R.1**) implies that $\Psi(z_k) = \Psi(\bar{z})$ holds for all $k \ge K$. Then $z_k$ is constant for $k \ge K$, hence has finite length.

On the other hand, let $\Psi(z_k) > \Psi(\bar{z})$, denote $\psi_k = \Psi(z_k) - \Psi(\bar{z})$. Owing to (**R.6**), (**R.7**) and Definition 2.1, the KL property of $\Psi$ means that there exist $\epsilon, \eta$ and a concave function $\varphi$, and

$$\mathcal{U} \overset{\text{def}}{=} \big\{u \in \mathbb{R}_s^n : \mathrm{dist}(u, \mathcal{C}_{z_k}) < \epsilon\big\} \bigcap \big[\Psi(\bar{z}) < \Psi(u) < \Psi(\bar{z}) + \eta\big], \tag{A.9}$$

sucht ath for all $z \in \mathcal{U}$,

$$\varphi'\big(\Psi(z) - \Psi(\bar{z})\big)\mathrm{dist}\big(0, \partial\Psi(z)\big) \ge 1. \tag{A.10}$$

Let $k_1 \in \mathbb{N}$ be such that $\Psi(z_k) < \Psi(\bar{z}) + \eta$ holds for all $k \ge k_1$. Owing to (**R.5**), there exists another $k_2 \in \mathbb{N}$ such that $\mathrm{dist}(z_k, \mathcal{C}_{z_k}) < \epsilon$ holds for all $k \ge k_2$. Let $K = \max\{k_1, k_2\}$, then $z_k \in \mathcal{U}$ holds for all $k \ge K$. Then from (A.10), we have for $k \ge K$

$$\varphi'(\psi_k)\mathrm{dist}\big(0, \partial\Psi(z_k)\big) \ge 1.$$

Since $\varphi$ is concave, $\varphi'$ is decreasing, and $\Psi(z_k)$ is decreasing, we have

$$\varphi(\psi_k) - \varphi(\psi_{k+1}) \ge \varphi'(\psi_k)\big(\Psi(z_k) - \Psi(z_{k+1})\big) \ge \frac{\Psi(z_k) - \Psi(z_{k+1})}{\mathrm{dist}(0, \partial\Psi(z_k))}.$$

From (**R.1**), since $\mathrm{dist}(0, \partial\Psi(z_k)) \le \|w_k\|$, then

$$\varphi(\psi_k) - \varphi(\psi_{k+1}) \ge \frac{\Psi(z_k) - \Psi(z_{k+1})}{\|w_k\|} \ge \frac{\Psi(z_k) - \Psi(z_{k+1})}{\sigma\sum_{j=k+1-s}^{k} \Delta_j}.$$

Moreover, $\Psi(z_k) - \Psi(z_{k+1}) \ge \delta\Delta_{k+1}^2$ from (**R.2**), therefore we get

$$\varphi(\psi_k) - \varphi(\psi_{k+1}) \ge \frac{\delta\Delta_{k+1}^2}{\sigma\sum_{j=k+1-s}^{k} \Delta_j},$$

which yields

$$\Delta_{k+1}^2 \le \left(\tfrac{\sigma}{\delta}(\varphi(\psi_k) - \varphi(\psi_{k+1}))\right) \textstyle\sum_{j=k+1-s}^{k} \Delta_j. \tag{A.11}$$

Taking the square root of both sides and applying Young's inequality with $\kappa > 0$, we further obtain

$$\begin{aligned} 2\Delta_{k+1} &\le \frac{1}{\kappa} \textstyle\sum_{j=k+1-s}^{k} \Delta_j + \dfrac{\kappa\sigma}{\delta}(\varphi(\psi_k) - \varphi(\psi_{k+1})) \\ (\kappa = s) \quad &\le \frac{1}{s} \textstyle\sum_{j=k+1-s}^{k} \Delta_j + \dfrac{s\sigma}{\delta}(\varphi(\psi_k) - \varphi(\psi_{k+1})). \end{aligned} \tag{A.12}$$

Summing up both sides over $k$, and since $x_0 = ... = x_{-s}$, we get

$$\ell \overset{\text{def}}{=} \textstyle\sum_{k\in\mathbb{N}} \Delta_k \le \Delta_1 + \dfrac{s\sigma}{\delta}\varphi(\psi_1) < +\infty,$$

which concludes the finite length property of $x_k$.

(ii) Then the convergence of the sequence follows from the fact that $\{x_k\}_{k\in\mathbb{N}}$ is a Cauchy sequence, hence convergent. Owing to Lemma A.6, there exists a critical point $x^\star \in \mathrm{crit}(\Phi)$ such that $\lim_{k\to\infty} x_k = x^\star$.

(iii) We now turn to proving local convergence to a global minimmizer. Note that if $x^\star$ is a global minimizer of $\Phi$, then $z^\star$ is a global minimizer of $\Psi$. Let $r > \rho > 0$ such that $\mathbb{B}_r(z^\star) \subset \mathcal{U}$ and $\eta < \delta(r-\rho)^2$. Suppose that the initial point $x_0$ is chosen such that following conditions hold,

$$\Psi(z^\star) \le \Psi(z_0) < \Psi(z^\star) + \eta \tag{A.13}$$

$$\|x_0 - x^\star\| + \ell(s-1) + 2\sqrt{\frac{\Psi(z_0) - \Psi(z^\star)}{\delta}} + \frac{\sigma}{\delta}\varphi(\psi_0) < \rho. \tag{A.14}$$

The descent property (**R.1**) of $\Psi$ together with (A.13) imply that for any $k \in \mathbb{N}$, $\Psi(z^\star) \le \Psi(z_{k+1}) \le \Psi(z_k) \le \Psi(z_0) < \Psi(z^\star) + \eta$, and

$$\|x_{k+1} - x_k\| \le \sqrt{\frac{\Psi(z_k) - \Psi(z_{k+1})}{\delta}} \le \sqrt{\frac{\Psi(z_k) - \Psi(z^\star)}{\delta}}. \tag{A.15}$$

Therefore, given any $k \in \mathbb{N}$, if we have $x_k \in \mathbb{B}_\rho(x^\star)$, then

$$\begin{aligned} \|x_{k+1} - x^\star\| &\le \|x_k - x^\star\| + \|x_{k+1} - x_k\| \le \|x_k - x^\star\| + \sqrt{\frac{\Psi(z_k) - \Psi(z^\star)}{\delta}} \\ &\le \rho + (r - \rho) = r, \end{aligned} \tag{A.16}$$

which means that $x_{k+1} \in \mathbb{B}_r(x^\star)$.

For any $k \in \mathbb{N}$, define the following partial sum

$$p_k = \textstyle\sum_{j=k+1-s}^{k-1} \sum_{i=1}^{j} \Delta_i.$$

Note that $p_k = 0$ for $k = 1$, and $\lim_{k\to\infty} p_k = \ell(s-1)$. Next we prove the following claims through induction: for $k \in \mathbb{N}$

$$x_k \in \mathbb{B}_\rho(x^\star) \tag{A.17}$$

$$\textstyle\sum_{j=1}^{k} \Delta_{j+1} + \Delta_{k+1} \le \Delta_1 + p_k + \frac{\sigma}{\delta}(\varphi(\psi_1) - \varphi(\psi_{k+1})). \tag{A.18}$$

From (A.15) we have

$$\|x_1 - x_0\| \le \sqrt{\frac{\Psi(z_0) - \Psi(z^\star)}{\delta}}. \tag{A.19}$$

Applying the triangle inequality we then obtain

$$\|x_1 - x^\star\| \le \|x_0 - x^\star\| + \|x_1 - x_0\| \le \|x_0 - x^\star\| + \sqrt{\frac{\Psi(z_0) - \Psi(z^\star)}{\delta}} < \rho,$$

which means $x_1 \in \mathbb{B}_\rho(x^\star)$. Now, taking $\kappa = 1$ in (A.12) yields, for any $k \in \mathbb{N}$,

$$2\Delta_{k+1} \le \sum\nolimits_{j=k+1-s}^{k} \Delta_j + \frac{\sigma}{\delta}(\varphi(\psi_k) - \varphi(\psi_{k+1})). \tag{A.20}$$

Let $k = 1$. Since $x_0 = \ldots = x_{-s}$, we have

$$2\Delta_2 \le \Delta_1 + \frac{\sigma}{\delta}(\varphi(\psi_1) - \varphi(\psi_2)).$$

Therefore, (A.17) and (A.18) hold for $k = 1$.
Now assume that they hold for some $k > 1$. Using the triangle inequality and (A.18),

$$\begin{aligned}
\|x_{k+1} - x^\star\| &\le \|x_0 - x^\star\| + \Delta_1 + \sum\nolimits_{j=1}^{k} \Delta_j \\
&\le \|x_0 - x^\star\| + 2\Delta_1 + p_k + \frac{\sigma}{\delta}(\varphi(\psi_1) - \varphi(\psi_{k+1})) \\
&\le \|x_0 - x^\star\| + 2\Delta_1 + (s-1)\ell + \frac{\sigma}{\delta}(\varphi(\psi_1) - \varphi(\psi_{k+1})) \\
\text{(A.19)} &\le \|x_0 - x^\star\| + 2\sqrt{\frac{\Psi(z_0) - \Psi(z^\star)}{\delta}} + (s-1)\ell + \frac{\sigma}{\delta}(\varphi(\psi_1) - \varphi(\psi_{k+1})).
\end{aligned}$$

As $\varphi(\psi) \ge 0$ and $\varphi'(\psi) > 0$ for $\psi \in ]0, \eta[$, and in view of (A.14), we arrive at

$$\|x_{k+1} - x^\star\| \le \|x_0 - x^\star\| + 2\sqrt{\frac{\Psi(z_0) - \Psi(z^\star)}{\delta}} + (s-1)\ell + \frac{\sigma}{\delta}\varphi(\psi_0) < \rho$$

whence we deduce that (A.17) holds at $k+1$. Now, taking (A.20) at $k+1$ gives

$$\begin{aligned}
2\Delta_{k+2} &\le \sum\nolimits_{j=k+2-s}^{k+1} \Delta_j + \frac{\sigma}{\delta}(\varphi(\psi_{k+1}) - \varphi(\psi_{k+2})) \\
&\le \Delta_{k+1} + \sum\nolimits_{j=k+2-s}^{k} \Delta_j + \frac{\sigma}{\delta}(\varphi(\psi_{k+1}) - \varphi(\psi_{(k+2)})).
\end{aligned} \tag{A.21}$$

Adding both sides of (A.21) and (A.18) we get

$$\begin{aligned}
\sum\nolimits_{j=1}^{k+1} \Delta_{j+1} + \Delta_{k+2} &\le \Delta_1 + p_k + \sum\nolimits_{j=k+2-s}^{k} \Delta_j + \frac{\sigma}{\delta}(\varphi(\psi_1) - \varphi(\psi_{k+2})) \\
&= \Delta_1 + p_{k+1} + \frac{\sigma}{\delta}(\varphi(\psi_1) - \varphi(\psi_{k+2})),
\end{aligned}$$

meaning that (A.18) holds at $k+1$. This concludes the induction proof.
In summary, the above result shows that if we start close enough from $x^\star$ (so that (A.13)-(A.14) hold), then the sequence $\{x_k\}_{k\in\mathbb{N}}$ will remain in the neighbourhood $\mathbb{B}_\rho(x^\star)$ and thus converges to a critical point $\bar{x}$ owing to Lemma A.6. Moreover, $\Psi(z_k) \to \Psi(\bar{z}) \ge \Psi(z^\star)$ by virtue of (**R.3**). Now we need to show that $\Psi(\bar{z}) = \Psi(z^\star)$. Suppose that $\Psi(\bar{z}) > \Psi(z^\star)$. As $\Psi$ has the KL property at $z^\star$, we have

$$\varphi'\big(\Psi(\bar{z}) - \Psi(z^\star)\big)\operatorname{dist}\big(0, \partial\Psi(\bar{z})\big) \ge 1.$$

But this is impossible since $\varphi'(s) > 0$ for $s \in ]0, \eta[$, and $\operatorname{dist}\big(0, \partial\Psi(\bar{z})\big) = 0$ as $\bar{z}$ is a critical point. Hence we have $\Psi(\bar{z}) = \Psi(z^\star)$, which means $\Phi(\bar{x}) = \Phi(x^\star)$, *i.e.* the cluster point $\bar{x}$ is actually a global minimizer. This concludes the proof.

□

# B Proof of Theorem 3.2

**Proof of Theorem 3.2.** Since $R \in \mathrm{PSF}_{x^\star}(\mathcal{M}_{x^\star})$ and $F$ is locally $C^2$ around $x^\star$, the smooth perturbation rule of partly smooth functions ensures that $\Phi \in \mathrm{PSF}_{x^\star}(\mathcal{M}_{x^\star})$ (Corollary 4.7 of [18]).

With conditions in Theorem 2.2 hold, we have there is a critical point $x^\star \in \mathrm{crit}(\Phi)$ such that $x_k \to x^\star$ and $\Phi(x_k) \to \Phi(x^\star)$ (proof of Lemma A.6).

The finite length property of $\{x_k\}_{k\in\mathbb{N}}$ gives $\Delta_k \to 0$, then from Lemma A.3, we have

$$\mathrm{dist}\big(0, \partial\Phi(x_{k+1})\big) \leq \|g_{k+1}\| \leq \big(\tfrac{1}{\underline{\gamma}} + L\big)\Delta_{k+1} + \textstyle\sum_{i\in I}\big(\tfrac{|a_{i,k}|}{\underline{\gamma}} + |b_{i,k}|\big)\Delta_{k-i}.$$

Altogether, this shows that the conditions of [21, Theorem 4.10] or [15, Proposition 10.12] are fulfilled on $R$ at $x^\star$ for $-\nabla F(x^\star)$, and the identification result follows. $\square$

# C Proof of Theorem 3.4

Before presenting the proofs, we need some extra result from partial smoothness, and also Riemannian geometry.

## C.1 Partial smoothness and Riemannian geometry

From the sharpness in Definition 3.1, Proposition 2.10 of [18] allows to prove the following fact.

**Fact C.1 (Local normal sharpness).** If $R \in \mathrm{PSF}_x(\mathcal{M})$, then all $x' \in \mathcal{M}$ near $x$ satisfy $\mathcal{T}_{\mathcal{M}}(x') = T_{x'}$. In particular, when $\mathcal{M}$ is affine or linear, then $T_{x'} = T_x$.

We now give expressions of the Riemannian gradient and Hessian (see Section C.2 for definitions) for the case of partly smooth functions relative to a $C^2$ submanifold. This is summarized in the following fact which follows by combining (C.2), (C.3), Definition 3.1, Fact C.1 and [13, Proposition 17] (or [25, Lemma 2.4]).

**Fact C.2.** If $R \in \mathrm{PSF}_x(\mathcal{M})$, then for any $x' \in \mathcal{M}$ near $x$

$$\nabla_{\mathcal{M}} R(x') = \mathrm{P}_{T_{x'}}(\partial R(x')),$$

and this does not depend on the smooth representation of $R$ on $\mathcal{M}$. In turn, for all $h \in T_{x'}$

$$\nabla^2_{\mathcal{M}} G(x')h = \mathrm{P}_{T_{x'}} \nabla^2 \widetilde{R}(x')h + \mathfrak{W}_{x'}\big(h, \mathrm{P}_{T_{x'}^{\perp}} \nabla \widetilde{R}(x')\big),$$

where $\widetilde{R}$ is a smooth extension (representative) of $R$ on $\mathcal{M}$, and $\mathfrak{W}_x(\cdot,\cdot) : T_x \times T_x^{\perp} \to T_x$ is the Weingarten map of $\mathcal{M}$ at $x$ (see Section C.2 below for definitions).

## C.2 Riemannian Geometry

Let $\mathcal{M}$ be a $C^2$-smooth embedded submanifold of $\mathbb{R}^n$ around a point $x$. With some abuse of terminology, we shall state $C^2$-manifold instead of $C^2$-smooth embedded submanifold of $\mathbb{R}^n$. The natural embedding of a submanifold $\mathcal{M}$ into $\mathbb{R}^n$ permits to define a Riemannian structure and to introduce geodesics on $\mathcal{M}$, and we simply say $\mathcal{M}$ is a Riemannian manifold. We denote respectively $\mathcal{T}_{\mathcal{M}}(x)$ and $\mathcal{N}_{\mathcal{M}}(x)$ the tangent and normal space of $\mathcal{M}$ at point near $x$ in $\mathcal{M}$.

**Exponential map** Geodesics generalize the concept of straight lines in $\mathbb{R}^n$, preserving the zero acceleration characteristic, to manifolds. Roughly speaking, a geodesic is locally the shortest path between two points on $\mathcal{M}$. We denote by $\mathfrak{g}(t;x,h)$ the value at $t \in \mathbb{R}$ of the geodesic starting at $\mathfrak{g}(0;x,h) = x \in \mathcal{M}$ with velocity $\dot{\mathfrak{g}}(t;x,h) = \frac{d\mathfrak{g}}{dt}(t;x,h) = h \in \mathcal{T}_{\mathcal{M}}(x)$ (which is uniquely defined). For every $h \in \mathcal{T}_{\mathcal{M}}(x)$, there exists an interval $I$ around 0 and a unique geodesic $\mathfrak{g}(t;x,h) : I \to \mathcal{M}$ such that $\mathfrak{g}(0;x,h) = x$ and $\dot{\mathfrak{g}}(0;x,h) = h$. The mapping

$$\mathrm{Exp}_x : \mathcal{T}_{\mathcal{M}}(x) \to \mathcal{M},\ \ h \mapsto \mathrm{Exp}_x(h) = \mathfrak{g}(1;x,h),$$

is called *Exponential map*. Given $x, x' \in \mathcal{M}$, the direction $h \in \mathcal{T}_{\mathcal{M}}(x)$ we are interested in is such that

$$\mathrm{Exp}_x(h) = x' = \mathfrak{g}(1;x,h).$$

**Parallel translation** Given two points $x, x' \in \mathcal{M}$, let $\mathcal{T}_{\mathcal{M}}(x), \mathcal{T}_{\mathcal{M}}(x')$ be their corresponding tangent spaces. Define

$$\tau : \mathcal{T}_{\mathcal{M}}(x) \to \mathcal{T}_{\mathcal{M}}(x'),$$

the parallel translation along the unique geodesic joining $x$ to $x'$, which is isomorphism and isometry w.r.t. the Riemannian metric.

**Riemannian gradient and Hessian** For a vector $v \in \mathcal{N}_{\mathcal{M}}(x)$, the Weingarten map of $\mathcal{M}$ at $x$ is the operator $\mathfrak{W}_x(\cdot, v) : \mathcal{T}_{\mathcal{M}}(x) \to \mathcal{T}_{\mathcal{M}}(x)$ defined by

$$\mathfrak{W}_x(\cdot, v) = -\mathrm{P}_{\mathcal{T}_{\mathcal{M}}(x)} \mathrm{d}V[h],$$

where $V$ is any local extension of $v$ to a normal vector field on $\mathcal{M}$. The definition is independent of the choice of the extension $V$, and $\mathfrak{W}_x(\cdot, v)$ is a symmetric linear operator which is closely tied to the second fundamental form of $\mathcal{M}$, see [12, Proposition II.2.1].

Let $G$ be a real-valued function which is $C^2$ along the $\mathcal{M}$ around $x$. The covariant gradient of $G$ at $x' \in \mathcal{M}$ is the vector $\nabla_{\mathcal{M}} G(x') \in \mathcal{T}_{\mathcal{M}}(x')$ defined by

$$\langle \nabla_{\mathcal{M}} G(x'),\, h \rangle = \frac{d}{dt} G\big(\mathrm{P}_{\mathcal{M}}(x' + th)\big)\big|_{t=0},\ \forall h \in \mathcal{T}_{\mathcal{M}}(x'),$$

where $\mathrm{P}_{\mathcal{M}}$ is the projection operator onto $\mathcal{M}$. The covariant Hessian of $G$ at $x'$ is the symmetric linear mapping $\nabla^2_{\mathcal{M}} G(x')$ from $\mathcal{T}_{\mathcal{M}}(x')$ to itself which is defined as

$$\langle \nabla^2_{\mathcal{M}} G(x')h,\, h \rangle = \frac{d^2}{dt^2} G\big(\mathrm{P}_{\mathcal{M}}(x' + th)\big)\big|_{t=0},\ \forall h \in \mathcal{T}_{\mathcal{M}}(x'). \tag{C.1}$$

This definition agrees with the usual definition using geodesics or connections [25]. Now assume that $\mathcal{M}$ is a Riemannian embedded submanifold of $\mathbb{R}^n$, and that a function $G$ has a $C^2$-smooth restriction on $\mathcal{M}$. This can be characterized by the existence of a $C^2$-smooth extension (representative) of $G$, *i.e.* a $C^2$-smooth function $\widetilde{G}$ on $\mathbb{R}^n$ such that $\widetilde{G}$ agrees with $G$ on $\mathcal{M}$. Thus, the Riemannian gradient $\nabla_{\mathcal{M}} G(x')$ is also given by

$$\nabla_{\mathcal{M}} G(x') = \mathrm{P}_{\mathcal{T}_{\mathcal{M}}(x')} \nabla \widetilde{G}(x'), \tag{C.2}$$

and $\forall h \in \mathcal{T}_{\mathcal{M}}(x')$, the Riemannian Hessian reads

$$\begin{aligned} \nabla^2_{\mathcal{M}} G(x')h &= \mathrm{P}_{\mathcal{T}_{\mathcal{M}}(x')} \mathrm{d}(\nabla_{\mathcal{M}} G)(x')[h] = \mathrm{P}_{\mathcal{T}_{\mathcal{M}}(x')} \mathrm{d}\big(x' \mapsto \mathrm{P}_{\mathcal{T}_{\mathcal{M}}(x')} \nabla_{\mathcal{M}} \widetilde{G}\big)[h] \\ &= \mathrm{P}_{\mathcal{T}_{\mathcal{M}}(x')} \nabla^2 \widetilde{G}(x')h + \mathfrak{W}_{x'}\big(h, \mathrm{P}_{\mathcal{N}_{\mathcal{M}}(x')} \nabla \widetilde{G}(x')\big), \end{aligned} \tag{C.3}$$

where the last equality comes from [1, Theorem 1]. When $\mathcal{M}$ is an affine or linear subspace of $\mathbb{R}^n$, then obviously $\mathcal{M} = x + \mathcal{T}_{\mathcal{M}}(x)$, and $\mathfrak{W}_{x'}(h, \mathrm{P}_{\mathcal{N}_{\mathcal{M}}(x')}\nabla\widetilde{G}(x')) = 0$, hence (C.3) reduces to

$$\nabla^2_{\mathcal{M}} G(x') = \mathrm{P}_{\mathcal{T}_{\mathcal{M}}(x')}\nabla^2\widetilde{G}(x')\mathrm{P}_{\mathcal{T}_{\mathcal{M}}(x')}.$$

See [17, 12] for more materials on differential and Riemannian manifolds.

The following lemmas summarize two key properties that we will need throughout.

**Lemma C.3.** *Let $x \in \mathcal{M}$, and $x_k$ a sequence converging to $x$ in $\mathcal{M}$. Denote $\tau_k : \mathcal{T}_{\mathcal{M}}(x) \to \mathcal{T}_{\mathcal{M}}(x_k)$ be the parallel translation along the unique geodesic joining $x$ to $x_k$. Then, for any bounded vector $u \in \mathbb{R}^n$, we have*

$$(\tau_k^{-1}\mathrm{P}_{\mathcal{T}_{\mathcal{M}}(x_k)} - \mathrm{P}_{\mathcal{T}_{\mathcal{M}}(x)})u = o(\|u\|).$$

**Proof.** See Lemma B.1 of [23]. □

**Lemma C.4.** *Let $x, x'$ be two close points in $\mathcal{M}$, denote $\tau : \mathcal{T}_{\mathcal{M}}(x) \to \mathcal{T}_{\mathcal{M}}(x')$ the parallel translation along the unique geodesic joining $x$ to $x'$. The Riemannian Taylor expansion of $\Phi \in C^2(\mathcal{M})$ around $x$ reads,*

$$\tau^{-1}\nabla_{\mathcal{M}}\Phi(x') = \nabla_{\mathcal{M}}\Phi(x) + \nabla^2_{\mathcal{M}}\Phi(x)\mathrm{P}_{\mathcal{T}_{\mathcal{M}}(x)}(x' - x) + o(\|x' - x\|).$$

**Proof.** See Lemma B.2 of [23]. □

## C.3 Proof of Theorem 3.4

The proof of Theorem 1 consists of several steps, first we prove that under the required setting, we can obtain (C.4), *i.e.* the linearized fixed-point iteration.

**Proposition C.5 (Locally linearized iteration).** *For Algorithm 1, suppose that conditions in Theorem 2.2 hold and the generated sequence $x_k$ converges to a critical point $x^\star \in \mathrm{crit}(\Phi)$ such that Theorem 3.2 and condition (3.2) and (3.3) hold. Then for all $k$ large enough, we have*

$$d_{k+1} = M d_k + o(\|d_k\|). \tag{C.4}$$

*The term $o(\cdot)$ vanishes if $R$ is polyhedral around $x^\star$ and $(\gamma_k, a_{i,k}, b_{i,k})$ are chosen constants.*

Define the iteration-dependent versions of the matrices in (3.1) and (3.4), *i.e.*

$$\begin{aligned}
&H_k \stackrel{\text{def}}{=} \gamma_k \mathrm{P}_{T_{x^\star}}\nabla^2 F(x^\star)\mathrm{P}_{T_{x^\star}},\ G_k \stackrel{\text{def}}{=} \mathrm{Id} - H_k,\ Q_k \stackrel{\text{def}}{=} \gamma_k\nabla^2_{\mathcal{M}_{x^\star}}\Phi(x^\star)\mathrm{P}_{T_{x^\star}} - H_k,\\
&M_{k,0} \stackrel{\text{def}}{=} (a_{k,0} - b_{k,0})P + (1 + b_{k,0})PG,\ M_{k,s} \stackrel{\text{def}}{=} -(a_{k,s-1} - b_{k,s-1})P - b_{k,s-1}PG,\\
&M_{k,i} \stackrel{\text{def}}{=} -\big((a_{k,i-1} - a_{k,i}) - (b_{k,i-1} - b_{k,i})\big)P - (b_{k,i-1} - b_{k,i})PG,\ i = 1, ..., s-1,\\
&M_k \stackrel{\text{def}}{=} \begin{bmatrix} M_{k,0} & M_{k,1} & \cdots & M_{k,s-1} & M_{k,s} \\ \mathrm{Id} & 0 & \cdots & 0 & 0 \\ 0 & \mathrm{Id} & \cdots & 0 & 0 \\ \vdots & \vdots & \ddots & \vdots & \vdots \\ 0 & 0 & \cdots & \mathrm{Id} & 0 \end{bmatrix}.
\end{aligned} \tag{C.5}$$

After the finite identification of $\mathcal{M}_{x^\star}$, we have $x_k \in \mathcal{M}_{x^\star}$ for $x_k$ close enough to $x^\star$. Let $T_{x_k}$ be their corresponding tangent spaces, and define $\tau_k : T_{x^\star} \to T_{x_k}$ the parallel translation along the unique geodesic joining from $x_k$ to $x^\star$.

Before proving Proposition C.5, we first establish the following intermediate result which provides useful estimates.

**Proposition C.6.** *Under the assumptions of Proposition C.5, we have*

$$\begin{aligned} &\|y_{a,k} - x^\star\| = O(\|d_k\|),\ \|y_{b,k} - x^\star\| = O(\|d_k\|),\ \|r_{k+1}\| = O(\|d_k\|), \\ &(\tau_{k+1}^{-1}\mathrm{P}_{T_{x_{k+1}}} - \mathrm{P}_{T_{x^\star}})\big(\nabla F(y_{b,k}) - \nabla F(x_{k+1})\big) = o(\|d_k\|). \end{aligned} \tag{C.6}$$

*and*

$$\|P(Q_k - Q)r_{k+1}\| = o(\|d_k\|),\ \|(M_k - M)d_k\| = o(\|d_k\|). \tag{C.7}$$

**Proof.** Since $|a_{i,k}| \le 1$, then

$$\begin{aligned} \|y_{a,k} - x^\star\| &= \|x_k + \textstyle\sum_{i\in I} a_{i,k}(x_{k-i} - x_{k-i-1}) - x^\star + \sum_{i\in I} a_{i,k}(x^\star - x^\star)\| \\ &\le \|x_k - x^\star\| + \textstyle\sum_{i\in I} a_{i,k}(\|x_{k-i} - x^\star\| + \|x_{k-i-1} - x^\star\|) \\ &\le 2\textstyle\sum_{i\in I} \|r_{k-i}\| \le 2\sqrt{s+1}\left\|\begin{matrix} r_k \\ \vdots \\ r_{k-s}\end{matrix}\right\| = 2\sqrt{s+1}\|d_k\|, \end{aligned} \tag{C.8}$$

hence we get the first and second estimates. From prox-regularity of $R$ at $x^\star$ for $-\nabla F(x^\star)$, invoking [30, Proposition 13.37], we have that there exists $\bar{r} > 0$ such that for all $\gamma_k \in ]0, \min(\bar{\gamma}, \bar{r})[$, there exists a neighbourhood $U$ of $x^\star - \gamma_k \nabla F(x^\star)$ on which $\mathrm{prox}_{\gamma_k R}$ is single-valued and $l$-Lipschitz continuous with $l = \bar{r}/(\bar{r} - \gamma_k)$. Since $\nabla F$ is continuous and $x_k \to x^\star$, we have $y_{a,k} - \gamma_k \nabla F(y_{b,k}) \to x^\star - \gamma_k \nabla F(x^\star)$. In turn, $y_{a,k} - \gamma_k \nabla F(y_{b,k}) \in U$ for all $k$ sufficiently large. Then in turn, we obtain

$$\begin{aligned} \|r_{k+1}\| &= \|\mathrm{prox}_{\gamma_k R}\big(y_{a,k} - \gamma_k \nabla F(y_{b,k})\big) - \mathrm{prox}_{\gamma_k R}\big(x^\star - \gamma_k \nabla F(x^\star)\big)\| \\ &\le l\|(y_{a,k} - x^\star) - \gamma_k\big(\nabla F(y_{b,k}) - \nabla F(x^\star)\big)\| \\ &\le l(\|y_{a,k} - x^\star\| + \gamma_k L\|y_{b,k} - x^\star\|) \\ &\le 2l\sqrt{s+1}(1 + \gamma_k L)\|d_k\| \le 4l\sqrt{s+1}\|d_k\|, \end{aligned} \tag{C.9}$$

which yields the third estimate. Combining Lemma C.3, (C.8) and (C.9), we get

$$\begin{aligned} (\tau_{k+1}^{-1}\mathrm{P}_{T_{x_{k+1}}} - \mathrm{P}_{T_{x^\star}})\big(\nabla F(y_{b,k}) - \nabla F(x_{k+1})\big) &= o(\|\nabla F(y_{b,k}) - \nabla F(x_{k+1})\|) \\ &= o(\|y_{b,k} - x^\star\|) + o(\|r_{k+1}\|) = o(\|d_k\|). \end{aligned}$$

Let's now turn to (C.7). First, define the function $\overline{R}(x) \stackrel{\mathrm{def}}{=} R(x) + \langle x, \nabla F(x^\star)\rangle$. From the smooth perturbation rule of partial smoothness [18, Corollary 4.7], $\overline{R} \in \mathrm{PSF}_{x^\star}(\mathcal{M}_{x^\star})$. Moreover, from Fact C.2 and normal sharpness, the Riemannian Hessian of $\overline{R}$ at $x^\star$ is such that, $\forall h \in T_{x^\star}$,

$$\begin{aligned} \gamma\nabla^2_{\mathcal{M}_{x^\star}}\overline{R}(x^\star)h &= \gamma\mathrm{P}_{T_{x^\star}}\nabla^2\widetilde{\overline{R}}(x^\star)h + \gamma\mathfrak{W}_{x^\star}\big(h, \mathrm{P}_{T_{x^\star}^\perp}\nabla\widetilde{\overline{R}}(x^\star)\big) \\ &= \gamma\mathrm{P}_{T_{x^\star}}\nabla^2\widetilde{R}(x^\star)h + \gamma\mathfrak{W}_{x^\star}\big(h, \mathrm{P}_{T_{x^\star}^\perp}\nabla\widetilde{\Phi}(x^\star)\big) \\ &= \gamma\nabla^2_{\mathcal{M}_{x^\star}}\Phi(x^\star)\mathrm{P}_{T_{x^\star}}h - Hh = Qh, \end{aligned}$$

where $\widetilde{\cdot}$ is the smooth representative of the corresponding function. We have

$$\begin{aligned} \lim_{k\to\infty}\frac{\|P(Q_k - Q)r_{k+1}\|}{\|r_{k+1}\|} &= \lim_{k\to\infty}\frac{\|P(\gamma_k - \gamma)\nabla^2_{\mathcal{M}_{x^\star}}\overline{R}(x^\star)\mathrm{P}_{T_{x^\star}}r_{k+1}\|}{\|r_{k+1}\|} \\ &\le \lim_{k\to\infty}|\gamma_k - \gamma|\|P\|\|\nabla^2_{\mathcal{M}_{x^\star}}\overline{R}(x^\star)\mathrm{P}_{T_{x^\star}}\| = 0, \end{aligned}$$

which entails $\|P(Q_k - Q)r_{k+1}\| = o(\|r_{k+1}\|) = o(\|d_k\|)$. Similarly, since $H$ is Lipschitz, we have

$$\lim_{k\to\infty} \frac{\|P(G_k - G)r_k\|}{\|r_k\|} = \lim_{k\to\infty} \frac{\|P(\gamma_k - \gamma)Hr_k\|}{\|r_k\|} \le \lim_{k\to\infty} |\gamma_k - \gamma|L\|P\| = 0. \tag{C.10}$$

Now, let's consider $(M_k - M)d_k$

$$M_k - M = \begin{bmatrix} M_{k,0} - M_0 & M_{k,1} - M_1 & \cdots & M_{k,s-1} - M_{s-1} & M_{k,s} - M_s \\ 0 & 0 & \cdots & 0 & 0 \\ 0 & 0 & \cdots & 0 & 0 \\ \vdots & \vdots & \ddots & \vdots & \vdots \\ 0 & 0 & \cdots & 0 & 0 \end{bmatrix}.$$

Take $(M_{k,0} - M_0)r_k$, we have

$$\begin{aligned} &(M_{k,0} - M_0)r_k \\ &= \big((a_{k,0} - b_{k,0})P + (1 + b_{k,0})PG_k\big)r_k - \big((a_0 - b_0)P + (1 + b_0)PG\big)r_k \\ &= \big((a_{k,0} - b_{k,0}) - (a_0 - b_0)\big)Pr_k + (1 + b_{k,0})P(G_k - G)r_k + (b_{k,0} - b_0)PGr_k. \end{aligned}$$

Since we assume that $a_{i,k} \to a_i, b_{i,k} \to b_i, i = 0, 1$ and $\gamma_k \to \gamma$, plus (C.10), it can be shown that

$$\begin{aligned} &\lim_{k\to\infty} \frac{\|(M_{k,0} - M_0)r_k\|}{\|r_k\|} \\ &\le \lim_{k\to\infty} |(a_{k,0} - b_{k,0}) - (a_0 - b_0)|\|P\| + |1 + b_{k,0}||\gamma_k - \gamma|L\|P\| + |b_{k,0} - b_0|\|P\|\|G\| = 0, \end{aligned}$$

that is $\|(M_{k,0} - M_0)r_k\| = o(\|r_k\|)$. Using the same arguments, we can show that

$$\|(M_{k,i} - M_i)r_{k-i}\| = o(\|r_{k-i}\|),\ i = 1, ..., s-1 \ \text{ and } \ \|(M_{k,s} - M_s)r_{k,s}\| = o(\|r_{k,s}\|).$$

Assemble them together, we obtain

$$\|(M_k - M)d_k\| = o(\|d_k\|),$$

which concludes the proof. □

**Proof of Proposition C.5.** From the update (1.7) and the condition for a critical point $x^\star$ of problem ($\mathcal{P}$), we have

$$\begin{aligned} y_{a,k} - x_{k+1} - \gamma_k\big(\nabla F(y_{b,k}) - \nabla F(x_{k+1})\big) &\in \gamma_k \partial\Phi(x_{k+1}) \\ 0 &\in \gamma_k \partial\Phi(x^\star). \end{aligned}$$

Projecting into $T_{x_{k+1}}$ and $T_{x^\star}$, respectively, and using Fact C.2, leads to

$$\begin{aligned} \gamma_k \tau_{k+1}^{-1} \nabla_{\mathcal{M}_{x^\star}} \Phi(x_{k+1}) &= \tau_{k+1}^{-1} \mathrm{P}_{T_{x_{k+1}}} \big(y_{a,k} - x_{k+1} - \gamma_k\big(\nabla F(y_{b,k}) - \nabla F(x_{k+1})\big)\big) \\ \gamma_k \nabla_{\mathcal{M}_{x^\star}} \Phi(x^\star) &= 0. \end{aligned}$$

Adding both identities, and subtracting $\tau_{k+1}^{-1} \mathrm{P}_{T_{x_{k+1}}} x^\star$ on both sides, we arrive at

$$\begin{aligned} &\tau_{k+1}^{-1} \mathrm{P}_{T_{x_{k+1}}} r_{k+1} + \gamma_k \big(\tau_{k+1}^{-1} \nabla_{\mathcal{M}_{x^\star}} \Phi(x_{k+1}) - \nabla_{\mathcal{M}_{x^\star}} \Phi(x^\star)\big) \\ &= \tau_{k+1}^{-1} \mathrm{P}_{T_{x_{k+1}}} (y_{a,k} - x^\star) - \gamma_k \tau_{k+1}^{-1} \mathrm{P}_{T_{x_{k+1}}} \big(\nabla F(y_{b,k}) - \nabla F(x_{k+1})\big). \end{aligned} \tag{C.11}$$

By virtue of Lemma C.3, we get

$$\tau_{k+1}^{-1}\mathrm{P}_{T_{x_{k+1}}} r_{k+1} = \mathrm{P}_{T_{x^\star}} r_{k+1} + (\tau_{k+1}^{-1}\mathrm{P}_{T_{x_{k+1}}} - \mathrm{P}_{T_{x^\star}}) r_{k+1} = \mathrm{P}_{T_{x^\star}} r_{k+1} + o(\|r_{k+1}\|).$$

Using [22, Lemma 5.1], we also have

$$r_{k+1} = \mathrm{P}_{T_{x^\star}} r_{k+1} + o(\|r_{k+1}\|),$$

and thus

$$\tau_{k+1}^{-1}\mathrm{P}_{T_{x_{k+1}}} r_{k+1} = r_{k+1} + o(\|r_{k+1}\|) = r_{k+1} + o(\|d_k\|), \tag{C.12}$$

where we also used (C.6). Similarly

$$\begin{aligned}
&\tau_{k+1}^{-1}\mathrm{P}_{T_{x_{k+1}}}(y_{a,k} - x^\star) \\
&= \mathrm{P}_{T_{x^\star}}(y_{a,k} - x^\star) + (\tau_{k+1}^{-1}\mathrm{P}_{T_{x_{k+1}}} - \mathrm{P}_{T_{x^\star}})(y_{a,k} - x^\star) \\
&= \mathrm{P}_{T_{x^\star}}(y_{a,k} - x^\star) + o(\|y_{a,k} - x^\star\|) = \mathrm{P}_{T_{x^\star}}(y_{a,k} - x^\star) + o(\|d_k\|) \\
&= \mathrm{P}_{T_{x^\star}}(x_k - x^\star) + \textstyle\sum_{i\in I} a_{i,k}\mathrm{P}_{T_{x^\star}}\big((x_{k-i} - x^\star) - (x_{k-i-1} - x^\star)\big) + o(\|d_k\|) \\
&= r_k + o(\|r_k\|) + \textstyle\sum_{i\in I} a_{i,k}\big(r_{k-i} - r_{k-i-1} + o(\|r_{k-i}\|) + o(\|r_{k-i-1}\|)\big) + o(\|d_k\|) \\
&= r_k + \textstyle\sum_{i\in I} a_{i,k}(r_{k-i} - r_{k-i-1}) + \sum_{i\in I\cup\{s\}} o(\|r_{k-i}\|) + o(\|d_k\|) \\
&= (y_{a,k} - x^\star) + o(\|d_k\|).
\end{aligned} \tag{C.13}$$

Moreover owing to Lemma C.4 and (C.6),

$$\begin{aligned}
\tau^{-1}\nabla_{\mathcal{M}_{x^\star}}\Phi(x_{k+1}) - \nabla_{\mathcal{M}_{x^\star}}\Phi(x^\star) &= \nabla^2_{\mathcal{M}_{x^\star}}\Phi(x^\star)\mathrm{P}_{T_{x^\star}} r_{k+1} + o(\|r_{k+1}\|) \\
&= \nabla^2_{\mathcal{M}_{x^\star}}\Phi(x^\star)\mathrm{P}_{T_{x^\star}} r_{k+1} + o(\|d_k\|).
\end{aligned} \tag{C.14}$$

Therefore, inserting (C.12), (C.13) and (C.14) into (C.11), we obtain

$$\begin{aligned}
&\big(\mathrm{Id} + \gamma_k\nabla^2_{\mathcal{M}_{x^\star}}\Phi(x^\star)\mathrm{P}_{T_{x^\star}}\big) r_{k+1} \\
&= (y_{a,k} - x^\star) - \gamma_k\tau_{k+1}^{-1}\mathrm{P}_{T_{x_{k+1}}}\big(\nabla F(y_{b,k}) - \nabla F(x_{k+1})\big) + o(\|d_k\|).
\end{aligned} \tag{C.15}$$

Owing to (C.6) and local $C^2$-smoothness of $F$, we have

$$\begin{aligned}
&\tau_{k+1}^{-1}\mathrm{P}_{T_{x_{k+1}}}\big(\nabla F(y_{b,k}) - \nabla F(x_{k+1})\big) \\
&= \mathrm{P}_{T_{x^\star}}\big(\nabla F(y_{b,k}) - \nabla F(x_{k+1})\big) + o(\|d_k\|) \\
&= \mathrm{P}_{T_{x^\star}}\big(\nabla F(y_{b,k}) - \nabla F(x^\star)\big) - \mathrm{P}_{T_{x^\star}}\big(\nabla F(x_{k+1}) - \nabla F(x^\star)\big) + o(\|d_k\|) \\
&= \mathrm{P}_{T_{x^\star}}\nabla^2 F(x^\star)(y_{b,k} - x^\star) + o(\|y_{b,k} - x^\star\|) - \mathrm{P}_{T_{x^\star}}\nabla^2 F(x^\star) r_{k+1} + o(\|r_{k+1}\|) + o(\|d_k\|) \\
&= \mathrm{P}_{T_{x^\star}}\nabla^2 F(x^\star)\mathrm{P}_{T_{x^\star}}(y_{b,k} - x^\star) - \mathrm{P}_{T_{x^\star}}\nabla^2 F(x^\star)\mathrm{P}_{T_{x^\star}}(x_{k+1} - x^\star) + o(\|d_k\|).
\end{aligned} \tag{C.16}$$

Injecting (C.16) in (C.15), we get

$$\begin{aligned}
&\big(\mathrm{Id} + \gamma_k\nabla^2_{\mathcal{M}_{x^\star}}\Phi(x^\star)\mathrm{P}_{T_{x^\star}} - \gamma_k\mathrm{P}_{T_{x^\star}}\nabla^2 F(x^\star)\mathrm{P}_{T_{x^\star}}\big) r_{k+1} \\
&= (\mathrm{Id} + Q_k) r_{k+1} = (y_{a,k} - x^\star) - H_k(y_{b,k} - x^\star) + o(\|d_k\|),
\end{aligned} \tag{C.17}$$

which can be further written as, recall that $H_k = \mathrm{Id} - G_k$,

$$
\begin{aligned}
&(\mathrm{Id} + Q_k) r_{k+1} \\
&= (\mathrm{Id} + Q) r_{k+1} + (Q_k - Q) r_{k+1} \\
&= (y_{a,k} - x^\star) - H_k (y_{b,k} - x^\star) + o(\|d_k\|) \\
&= r_k + \textstyle\sum_{i\in I} a_{i,k}(r_{k-i} - r_{k-i-1}) - H_k\big(r_k + \textstyle\sum_{i\in I} b_{i,k}(r_{k-i} - r_{k-i-1})\big) + o(\|d_k\|) \\
&= (1 + a_{k,0}) r_k - \textstyle\sum_{i=1}^{s-1} (a_{k,i-1} - a_{k,i}) r_{k-i} - a_{k,s-1} r_{k-s} \\
&\quad - H_k\big((1 + b_{k,0}) r_k - \textstyle\sum_{i=1}^{s-1} (b_{k,i-1} - b_{k,i}) r_{k-i} - b_{k,s-1} r_{k-s}\big) + o(\|d_k\|) \\
&= (1 + a_{k,0}) r_k - \textstyle\sum_{i=1}^{s-1} (a_{k,i-1} - a_{k,i}) r_{k-i} - a_{k,s-1} r_{k-s} \\
&\quad - (1 + b_{k,0}) H_k r_k + H_k \textstyle\sum_{i=1}^{s-1} (b_{k,i-1} - b_{k,i}) r_{k-i} + H_k b_{k,s-1} r_{k-s} + o(\|d_k\|) \\
&= \big((1 + a_{k,0})\mathrm{Id} - (1 + b_{k,0}) H_k\big) r_k - (a_{k,s-1}\mathrm{Id} - b_{k,s-1} H_k) r_{k-s} \\
&\quad - \textstyle\sum_{i=1}^{s-1} \big((a_{k,i-1} - a_{k,i})\mathrm{Id} - (b_{k,i-1} - b_{k,i}) H_k\big) r_{k-i} + o(\|d_k\|) \\
&= \big((a_{k,0} - b_{k,0})\mathrm{Id} + (1 + b_{k,0}) G_k\big) r_k - \big((a_{k,s-1} - b_{k,s-1})\mathrm{Id} + b_{k,s-1} G_k\big) r_{k-s} \\
&\quad - \textstyle\sum_{i=1}^{s-1} \big((a_{k,i-1} - a_{k,i})\mathrm{Id} - (b_{k,i-1} - b_{k,i})\mathrm{Id} + (b_{k,i-1} - b_{k,i}) G_k\big) r_{k-i} + o(\|d_k\|).
\end{aligned}
$$

Inverting $\mathrm{Id} + Q$ (which is possible thanks to assumption (3.2)), we obtain

$$
\begin{aligned}
&r_{k+1} + P(Q_k - Q) r_{k+1} \\
&= \big((a_{k,0} - b_{k,0}) P + (1 + b_{k,0}) P G_k\big) r_k - \big((a_{k,s-1} - b_{k,s-1}) P + b_{k,s-1} P G_k\big) r_{k-s} \\
&\quad - \textstyle\sum_{i=1}^{s-1} \big((a_{k,i-1} - a_{k,i}) P - (b_{k,i-1} - b_{k,i}) P + (b_{k,i-1} - b_{k,i}) P G_k\big) r_{k-i} + o(\|d_k\|) \\
&= M_{k,0} r_k + M_{k,s} r_{k-s} + \textstyle\sum_{i=1}^{s-1} M_{k,i} r_{k-i} + o(\|d_k\|).
\end{aligned}
$$

Using the estimates (C.7), we get

$$ d_{k+1} = (M + (M_k - M)) d_k + o(\|d_k\|) = M d_k + o(\|d_k\|). $$

□

With the above result, we are able to prove the claim (3.6), hence Theorem 3.4.

**Proof of Theorem 3.4.** Since $\rho(M) < 1$, then we have $M$ is convergent with $\lim_{k\to\infty} M^k = 0$. Define $\psi_k = o(d_k)$, suppose after $K > 0$ iterations, (C.4) holds, then for $k \geq K$

$$ d_{k+1} = M^{k+1-K} d_K + \textstyle\sum_{j=K}^{k} M^{k-j} \psi_j \tag{C.18} $$

Since the spectral radius $\rho(M) < 1$, then from the spectral radius formula, given any $\rho \in ]\rho(M), 1[$, there exists a constant $C$ such that, for any $k \in \mathbb{N}$

$$ \|M^k\| \leq \|M\|^k \leq C\rho^k. $$

Therefore, from (C.18), we get

$$
\begin{aligned}
\|d_{k+1}\| &\leq \|M^{k+1-K} d_K + \textstyle\sum_{j=K}^{k} M^{k-j} \psi_j\| \\
&\leq \|M\|^{k+1-K} \|d_K\| + \textstyle\sum_{j=K}^{k} \|M\|^{k-j} \|\psi_j\| \\
&\leq C\rho^{k+1-K} \|d_K\| + C \textstyle\sum_{j=K}^{k} \rho^{k-j} \|\psi_j\|.
\end{aligned}
$$

Together with the fact that $\psi_j = o(\|d_j\|)$ leads to the claimed result. See also the result of [29, Section 2.1.2, Theorem 1]. □

## References

[1] P-A. Absil, R. Mahony, and J. Trumpf. An extrinsic look at the Riemannian Hessian. In *Geometric Science of Information*, pages 361–368. Springer, 2013.

[2] F. Alvarez. On the minimizing property of a second order dissipative system in Hilbert spaces. *SIAM Journal on Control and Optimization*, 38(4):1102–1119, 2000.

[3] F. Alvarez and H. Attouch. An inertial proximal method for maximal monotone operators via discretization of a nonlinear oscillator with damping. *Set-Valued Analysis*, 9(1-2):3–11, 2001.

[4] H. Attouch, J. Bolte, and B. F. Svaiter. Convergence of descent methods for semi-algebraic and tame problems: proximal algorithms, Forward–Backward splitting, and regularized Gauss–Seidel methods. *Mathematical Programming*, 137(1-2):91–129, 2013.

[5] H. Attouch, J. Peypouquet, and P. Redont. A dynamical approach to an inertial Forward–Backward algorithm for convex minimization. *SIAM J. Optim.*, 24(1):232–256, 2014.

[6] A. Beck and M. Teboulle. A fast iterative shrinkage-thresholding algorithm for linear inverse problems. *SIAM Journal on Imaging Sciences*, 2(1):183–202, 2009.

[7] J. Bolte, A. Daniilidis, and A. Lewis. The Łojasiewicz inequality for nonsmooth subanalytic functions with applications to subgradient dynamical systems. *SIAM Journal on Optimization*, 17(4):1205–1223, 2007.

[8] J. Bolte, A. Daniilidis, O. Ley, and L. Mazet. Characterizations of Lojasiewicz inequalities: subgradient flows, talweg, convexity. *Transactions of the American Mathematical Society*, 362(6):3319–3363, 2010.

[9] R. I. Boţ, E. R. Csetnek, and S. C. László. An inertial Forward–Backward algorithm for the minimization of the sum of two nonconvex functions. *EURO Journal on Computational Optimization*, pages 1–23, 2014.

[10] E. J. Candès, X. Li, Y. Ma, and J. Wright. Robust principal component analysis? *Journal of the ACM (JACM)*, 58(3):11, 2011.

[11] A. Chambolle and C. Dossal. On the convergence of the iterates of the "Fast Iterative Shrinkage/Thresholding Algorithm". *Journal of Optimization Theory and Applications*, pages 1–15, 2015.

[12] I. Chavel. *Riemannian geometry: a modern introduction*, volume 98. Cambridge University Press, 2006.

[13] A. Daniilidis, W. Hare, and J. Malick. Geometrical interpretation of the predictor-corrector type algorithms in structured optimization problems. *Optimization: A Journal of Mathematical Programming & Operations Research*, 55(5-6):482–503, 2009.

[14] D. L. Donoho, M. Elad, and V. N. Temlyakov. Stable recovery of sparse overcomplete representations in the presence of noise. *IEEE Trans. Inform. Theory*, 52(1):6–18, 2006.

[15] D. Drusvyatskiy and A. S. Lewis. Optimality, identifiability, and sensitivity. *Mathematical Programming*, pages 1–32, 2013.

[16] H. Y. Le. Generalized subdifferentials of the rank function. *Optimization Letters*, 7(4):731–743, 2013.

[17] J. M. Lee. *Smooth manifolds*. Springer, 2003.

[18] A. S. Lewis. Active sets, nonsmoothness, and sensitivity. *SIAM J. on Optimization*, 13(3):702–725, 2003.

[19] A. S. Lewis and J. Malick. Alternating projections on manifolds. *Mathematics of Operations Research*, 33(1):216–234, 2008.

[20] A. S. Lewis and H. S. Sendov. Twice differentiable spectral functions. *SIAM Journal on Matrix Analysis and Applications*, 23(2):368–386, 2001.

[21] A. S. Lewis and S. Zhang. Partial smoothness, tilt stability, and generalized Hessians. *SIAM Journal on Optimization*, 23(1):74–94, 2013.

[22] J. Liang, J. Fadili, and G. Peyré. Local linear convergence of Forward–Backward under partial smoothness. In *Advances in Neural Information Processing Systems*, pages 1970–1978, 2014.

[23] J. Liang, M. J. Fadili, and G. Peyré. Activity identification and local linear convergence of Forward–Backward-type methods. arXiv:1503.03703, 2015.

[24] D. A. Lorenz and T. Pock. An inertial Forward–Backward algorithm for monotone inclusions. *Journal of Mathematical Imaging and Vision*, 51(2):311–325, 2014.

[25] S. A. Miller and J. Malick. Newton methods for nonsmooth convex minimization: connections among-Lagrangian, Riemannian Newton and SQP methods. *Math. programming*, 104(2-3):609–633, 2005.

[26] A. Moudafi and M. Oliny. Convergence of a splitting inertial proximal method for monotone operators. *Journal of Computational and Applied Mathematics*, 155(2):447–454, 2003.

[27] P. Ochs, Y. Chen, T. Brox, and T. Pock. iPiano: inertial proximal algorithm for nonconvex optimization. *SIAM Journal on Imaging Sciences*, 7(2):1388–1419, 2014.

[28] B. T. Polyak. Some methods of speeding up the convergence of iteration methods. *USSR Computational Mathematics and Mathematical Physics*, 4(5):1–17, 1964.

[29] B. T. Polyak. *Introduction to optimization*. Optimization Software, 1987.

[30] R. T. Rockafellar and R. Wets. *Variational analysis*, volume 317. Springer Verlag, 1998.

[31] L. van den Dries. *Tame topology and o-minimal structures*, volume 248 of *Mathematrical Society Lecture Notes*. Cambridge Univiversity Press, New York, 1998.